\documentclass[12pt, a4paper]{article}
\usepackage{amsmath,amssymb,bm,graphicx,latexsym,graphics,epsfig,subfigure,color}
\usepackage{amsthm}
\usepackage{mathrsfs}
\usepackage{amssymb,amsmath,amsfonts,longtable,color}
\usepackage{amsbsy, amsmath, amssymb, multirow, epstopdf,epsf}
\usepackage{setspace}
\usepackage{amsthm}
\usepackage{graphicx,float,wasysym}
\usepackage{color}
\usepackage{enumerate} 
\usepackage{authblk} 

\usepackage{graphicx}			
\usepackage{mathrsfs}
\usepackage{amsfonts}
\usepackage{amsmath,amssymb,mathrsfs,bm,latexsym,graphicx,color}		
\usepackage{amssymb}
\usepackage[dvipsnames]{xcolor}
\usepackage{natbib} \bibpunct{(}{)}{;}{a}{,}{,} 
\usepackage[normalem]{ulem}
\usepackage[colorlinks=true,urlcolor=blue,citecolor=blue,linkcolor=blue,bookmarks=true]{hyperref}
\usepackage{tikz}
\usepackage{circuitikz}		
\usepackage{multicol}
\usepackage{multirow}

\usepackage[T1]{fontenc}
\usepackage[utf8]{inputenc}
\usepackage{indentfirst}
\usepackage{algorithm} 
\usepackage{algorithmic} 
\usepackage{amsmath}
\usepackage{graphicx}
\usepackage{algorithm} 
\usepackage{algorithmic} 
\usepackage{multirow} 
\usepackage{amsmath} 
\usepackage{xcolor}
\usepackage[figuresright]{rotating}

\usepackage{float}
\usepackage{afterpage}

\usepackage{graphicx,psfrag,epsfig}

\allowdisplaybreaks

\usepackage{tikz,xcolor,hyperref}

\definecolor{lime}{HTML}{A6CE39}
\DeclareRobustCommand{\orcidicon}{
	\begin{tikzpicture}
		\draw[lime, fill=lime] (0,0)
		circle[radius=0.16]
		node[white]{{\fontfamily{qag}\selectfont \tiny \.{I}D}}; 
	\end{tikzpicture}
	\hspace{-2mm}
}
\foreach \x in {A, ..., Z}{%
	\expandafter\xdef\csname orcid\x\endcsname{\noexpand\href{https://orcid.org/\csname orcidauthor\x\endcsname}{\noexpand\orcidicon}}
}
\usepackage{marvosym}

\title{\bf Statistical Analysis of Chen Distribution Under Improved Adaptive Type-II Progressive Censoring
}

\author{Li Zhang\Letter\thanks{\Letter The corresponding author. Email address: 2021212088@nwnu.edu.cn}\hspace{-0.8mm}\orcidB{}
	\\ 
{College of Mathematics and Statistics }\\{Northwest Normal University, Lanzhou 730070, China}\\ \medskip
}

\begin{document}
	\maketitle	
	\begin{abstract}


		This paper takes into account the estimation for the two unknown parameters of the Chen distribution with bathtub-shape hazard rate function under the improved adaptive Type-II progressive censored data. Maximum likelihood estimation for two parameters are proposed and the approximate confidence intervals are established using the asymptotic normality. Bayesian estimation are obtained under the symmetric and asymmetric loss function, during which the importance sampling and Metropolis-Hastings algorithm are proposed. Finally, the performance of various estimation methods is evaluated by Monte Carlo simulation experiments, and the proposed estimation method is illustrated through the analysis of a real data set.
		

	\end{abstract}
	\section{Introduction}

In reliability research and lifetime test experiments, most of the experiments take a long time to terminate, but considering the cost and time of the experiment, the failure time of all individuals cannot be observed, only the exact time of failure of a few experimental individuals can be observed, and the factors that may lead to individual failures need to be considered, the following research on censored data is reasonable and necessary. The classical censoring schemes includes type-I and type-II censoring, which are the most basic censoring schemes. They can be described as: suppose there are $n$ independent identical units placed in a lifetime experiment, type-I censoring requires that the experiment be terminated at a prefixed time point $T$, while type-II censoring requires that he experiment is terminated when the predetermined number of failed individuals $m<n$ is observed. The hybrid censoring scheme is a mixture of type-I censoring and type-II censoring. \cite{Epstein1954} proposed a type-I hybrid censoring scheme for the first time. \cite{Bhattacharyya1987} derived the exact distribution of the maximum likelihood estimator of the mean of an exponential distribution and an exact lower confidence bound for the mean based on a hybrid censored sample. \cite{childs2003} propose a hybrid censoring scheme which guarantees at least a fixed number of failures in a life testing experiment and the exacted distribution of maximum likelihood estimate with exponential distribution as well as exact lower confidence bound for the mean ia studied. \cite{Kundu2007} presents the statistical inference on Weibull parameters when the data are hybrid censored. The disadvantage of the above three censoring schemes is that these cells are not removed at all time points except at the termination of the experiment. To solve this problem, a type-II progressive censoring scheme is used, which is a scheme that combines type-II progressive censoring and type-I censoring. Refer to the type-II progressive censoring scheme proposed by \cite{Balakrishnan2000}. The type-II progressive censoring scheme can be described as: consider $n$ identical units placed in a lifetime test experiment, let $X_1, X_2, \cdots ,X_n$ be the corresponding failure time and $X_1< X_2< \cdots <X_n$, let $m$ is the predetermined number of failures, when the first failure has occurred, record the failure time as $X_{1:m:n}$, at this time there are $R_1$ units from the remaining $n-1$ units randomly removed. Similarly, when the second failure time $X_{2:m:n}$ is observed, $R_2$ units are randomly removed from the remaining $n-2-R_1$ units, and so on, at the $m$-$th$ failure time $X_{m:m:n}$, all remaining $n-m-R_1-R_2-\cdots-R_{m-1}$ units are removed. In the progressive censoring scheme, $R_1,R_2,\cdots ,R_m$ are predetermined prior to the study and are not changed during the experiment. \cite{Childs2008ExactLI} proposed a type-II progressive hybrid censoring scheme. If $X_{m:m:n}>T$, the experiment is terminated at $X_{m:m:n}$, otherwise, the experiment terminated at $T$. \cite{dey2014} takes into account the estimation for the unknown parameter of the Rayleigh distribution under type-II progressive censoring scheme. \cite{Sanjeev2015} considered point and interval estimation for the Maxwell distribution of type-I progressive hybrid censored data.\cite{Alma2016} discussed the analysis of progressive type-II progressive hybrid censored data when the lifetime distribution of the individual item is the normal and extreme value distributions. \cite{Cramer2016} considered the exponential distribution under this censoring scheme. \cite{Zhangandgui2019} studied the inference of reliability of generalized Rayleigh distribution based on the progressively type-II censored data. \cite{Zhanggui2020} develop a goodness of fit test process for Pareto distribution based on progressive type-II censoring scheme.  

However, a problem with this type of censoring scheme is that the number of effective units is random. If the effective sample obtained is very small or close to $0$, it will make the statistical inference process infeasible or the experimental efficiency is very low. Therefore, to avoid this problem, \cite{Ng2009} proposed for the first time a new scheme called adaptive type-II progressive censoring (AT-II PCS) which is a mixture of type-I censoring and type-II progressive censoring. Under this censoring scheme, the number of failed units $m$ to be observed and the progressive censoring scheme $R_1,R_2,\cdots,R_m$ are given in advance, but in the process of the experiment, some $R_i$ values may change according to the situation, so this censoring scheme can try to strike a balance between the total experiment time, the number of failed units, and the validity of the statistical analysis. The detailed description is as follows: if $X_{m:m:n}<T$, the experiment ends at point $X_{m:m:n}$, and the remaining $R_m$ units are all removed. Otherwise, if before time point $T$, the number of failed individuals is $j (0< j <m)$ and $X_{j:m:n}<T<X_{j+1:m:n}$. Let $R_{j+1}=\cdots=R_{m-1}=0$, The test continued until the failure of the m-th unit was observed and $R_m=n-m-\sum_{i=1}^{j} R_i \quad(i=1,\cdots,j)$, the experiment ends at point $X_{m:m:n}$. Here $T$, $m$, $R_1$, $R_2$, $\cdots$, $R_m$ $(R_1+R_2+\cdots+R_m+m=n)$ are all preset before the experiment. The value of time $T$ plays a very important role in determining the value of $R_i$. When $T\to\infty$, it is obvious that time is not the main concern of the experimenter. So the censoring scheme will degenerate into a general type-II progressive censoring, if $T\to0$, the censoring scheme degenerates into the traditional type-II censoring scheme. Some general statistical properties of an adaptive type-II progressive censoring scheme are investigated by \cite{Ye2014}. \cite{Nassar2017} describes the estimation of the inverse Weibull parameters under the AT-II PCS. \cite{MohieElDin2017} studied the estimation of generalized exponential distribution based on AT-II PCS. \cite{Hanieh2020} discussed the problem of estimating parameters of the inverted exponentiated Rayleigh distribution based on AT-II PCS. \cite{Rakhi2021} studied the estimation of parameters for a two parameter Kumaraswamy-exponential distribution based on AT-II PCS. An adaptive type-II progressive hybrid censored sampling with random removals is considered by \cite{Ahmed2022}. \cite{qigui2022} studied the statistical inference of Gompertz distribution based on AT-II PCS. 
In AT-II PCS, the time of the experiment is not an important consideration. We focus on being able to observe enough $m$ failure units. If the test unit is a high-reliability product, the duration of the experiment will be very long, while AT-II PCS does not guarantee a satisfactory, appropriate experimental total test time. But in many practical situations, the time of the experiment is inevitably an important factor to consider.

 To remedy this deficiency, \cite{YAN202138} proposed an improved censoring scheme called improved adaptive Type-II progressive censoring scheme(IAT-II PCS). One of its advantages is that it is guaranteed to end within the time specified in the experiment. A detailed description of IAT-II PCS is as follows: suppose there are $n$ units in the experiment. Assuming that they are independent and identically distributed, the failure unit $m$ to be observed and the censoring scheme $R=(R_1,R_2,...,R_m)$ are per-specified before the experiment. In addition, two time thresholds $T_1$ and $T_2$ $(T_1,T_2>0)$ are per-specified with $T_1<{T_2}$, where time $T_1$ is a warning about the test time, and $T_2$ is the maximum time allowed for the experiment, which means that when the experiment is carried out to time $T_1$, the experiment needs to be accelerated. In order to ensure that as many failed individuals as possible can be observed when the experiment is terminated at $T_2$, no experimental units will be censored when the experiment is accelerated. Of course, under this censoring scheme, when the experiment is carried out for enough time, the experiment is allowed to observe less than $m$ failure individuals. The experiment has the following three cases.  
 \begin{figure}[htbp]
 	\begin{flushleft}
 	\includegraphics[width=1.2\linewidth]{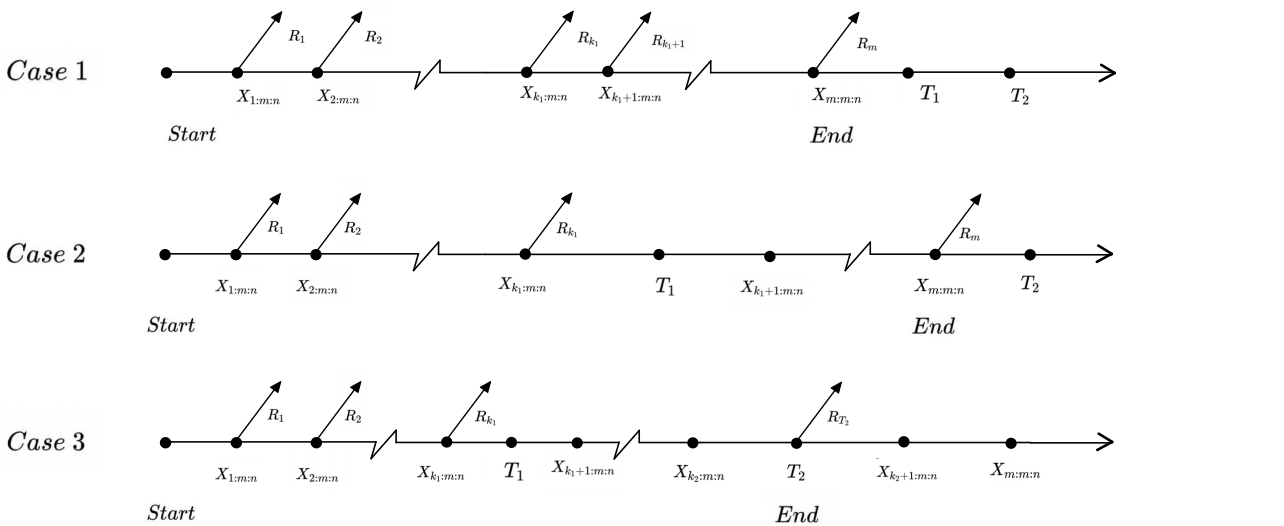}
 	\end{flushleft}
 	\caption{{\small Schematic representation of IAT-II PCS}}
 	\label{fig:case}
 \end{figure}
\begin{eqnarray}\label{case}
&Case1:&X_{1:m:n} , X_{2:m:n},\dots ,X_{m:m:n}, \quad if\  X_{m:m:n}<T_1<T_2, \\
	&Case2:&X_{1:m:n} ,\dots ,X_{k_1:m:n}\dots , X_{m:m:n}, \quad if\ X_{k_1:m:n}<T_1<X_{m:m:n}<T_2,\\
	&Case3:&X_{1:m:n} , \dots ,X_{k_1:m:n}\dots , X_{k_2:m:n}, \quad if\ X_{k_1:m:n}<T_1<X_{k_2:m:n}<T_2,	
\end{eqnarray}

The censoring scheme and the end position of the experiment in these three cases are as follows

\begin{eqnarray}\label{caseR}
	&Case1:&R=(R_1,R_2,...,R_m),R_m=0,\\ 
	&Case2:&R=(R_1,R_2,...,R_{k_1},0,...,0,R_m), R_m=n-m-\sum_{i=1}^{k_1} R_i,\\
	&Case3:&R=(R_1,R_2,...,R_{k_1},0,...,R_{k_2},R_{T_2}),R_{k_2}=0,R_{T_2}=n-k_2-\sum_{i=1}^{k_1} R_i,
\end{eqnarray}

In real life, there are many hazard rate functions used to model real life data, the most popular failure rate functions are constant, increasing or decreasing failure rate functions. For example Weibull, gamma and exponentiated exponential, among others. But there are other different forms of the hazard function, such as unimodal, bathtub-shaped, increasing-decreasing-increasing, etc. These distribution models described above also do not fit reasonably non-monotonic hazard rates, especially bathtub shape hazard rates that are common in reliability and other fields. An example of a bathtub-shaped failure rate is that the failure rate of newly installed electronic equipment is relatively high, failures often occur, production performance drops to a minimum level for a short period of time, and then the equipment will stabilize after a few days or months of use, remaining at this level for several years, and then the equipment underwent various forms of overhaul or technical transformation to restore production performance. \cite{YAN202138} used the Burr-XII distribution in the IAT-II PC proposed in 2021, and its hazard rate function could not fit the data of the bathtub curve. In recent years, a number of probability distributions have been introduced to model the risk rate of having a bathtub shape, which can be found in the literature. In this paper, we mainly focus on a two-parameter bathtub distribution proposed by \cite{CHEN2000155}. the probability density function $(PDF)$ and cumulative distribution function $(CDF)$ of $X$ are given by
\begin{eqnarray}\label{chenfF}
&f(x;\alpha ,\beta )&=\alpha \beta x^{\beta -1}\exp\left[ \alpha (1-e^{x^\beta })+x^\beta\right] \quad x>0,\alpha, \beta >0,\\
&F(x;\alpha ,\beta )&=1-\exp\left[ \alpha (1-e^{x^\beta })\right]  \quad x>0,\alpha, \beta >0.
\end{eqnarray}

The reliability and the hazard rate functions hazard function of $X$ are as follows
\begin{eqnarray}\label{chenSH}
S(x;\alpha ,\beta )&=&\exp\left[ \alpha (1-e^{x^\beta })\right] \quad x>0,\alpha, \beta >0,  \\
H(x;\alpha ,\beta )&=&\alpha \beta x^{\beta -1}e^{x^\beta } \quad x>0,\alpha, \beta >0,
\end{eqnarray}
where $\alpha$ and $\beta$ are the unknown scale and shape parameters, respectively, and when $\beta<1$, the hazard rate function is bathtub-shape.\\
\begin{figure*}
	\centering
	\includegraphics[width=1.15\linewidth]{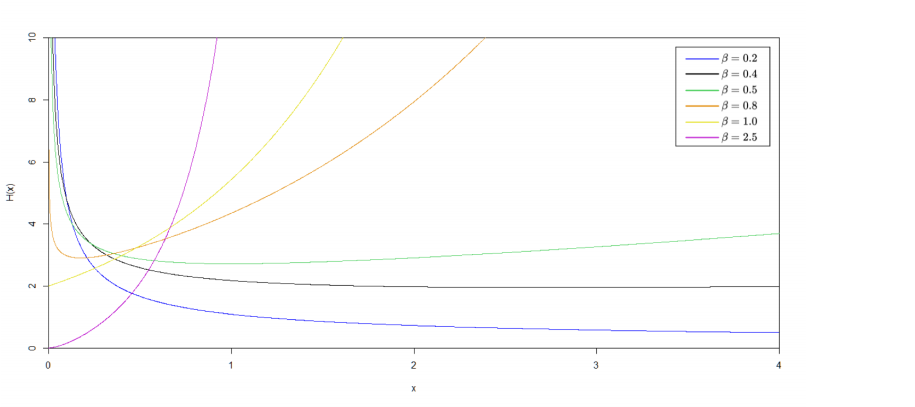}
	\caption{{\small Failure rate functions with $\alpha=2$.}}
	\label{fig:h2}
\end{figure*}
\begin{figure*}
	\centering
	\includegraphics[width=1.15\linewidth]{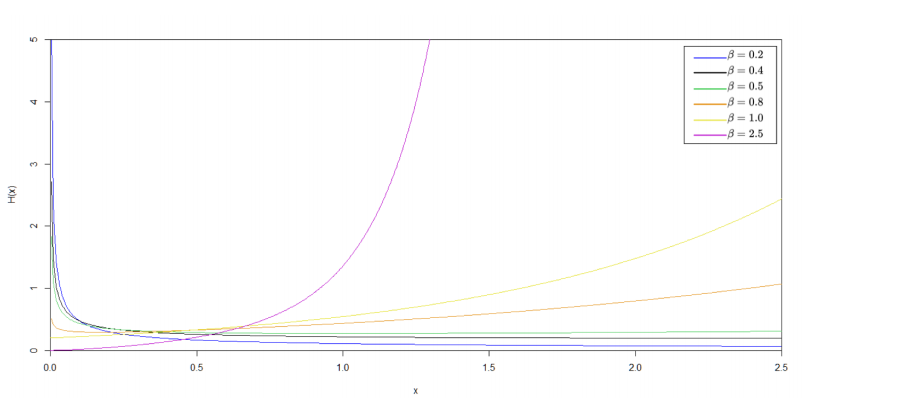}
	\caption{{\small Failure rate functions with $\alpha=0.2$.}}
	\label{fig:h0.2}
\end{figure*}

The rest of this paper is organized as follows: In Section $2$, maximum likelihood estimate (MLE) on unknown parameters is derived and the approximate confidence intervals is proposed. In Section $3$, Bayesian estimation is performed with gamma prior, and Importance sampling and Metropolis-Hastings algorithm are used. Next, based on the Monte Carlo method, a large number of simulations are carried out for parameter estimations and interval estimations in Section $4$. A real data set has been used to evaluated the estimation methods in Section $2$ and $3$. Finally, in Section $6$, we put forward some concluding conclusions.

\section{Model description and notation}
Let $X=X_{1:m:n}^R,X_{2:m:n}^R,\dots,X_{k_1:m:n}^R,\dots,X_{k_2:m:n}^R,\dots,X_{m:m:n}^R$ is an IAT-II PC sample from a lifetime test of size $m$ from a sample of size $n$ , where lifetimes have Chen distribution with pdf,cdf as given by \eqref{chenfF}. With predetermined number of removal of units from experiment, say $R=(R_1,R_2,\dots,R_{k_1},\dots,R_{k_2},\dots,R_m)$. The parameter $x_{i:m:n}^R$ (simplified as $x_i$ in later equation, $i=1,2,\dots,m$) is used to represent the observed values of IAT-II PC sample. on this basis, the corresponding likelihood function is given by 

%
\begin{eqnarray}\label{L_IAT- II PC}
	L=C\prod_{i=1}^{D_2}f(x_i)\prod_{i=1}^{D_1}\left[\left( 1-F(x_i)\right) \right]^{R_i}\left( 1-F(x_B)\right) ^{B},
\end{eqnarray}
where $D_2$, $D_1$, $B$, $C$ are shown in Table \eqref{xishu}

\begin{table}[H]
	\centering
	\setlength\tabcolsep{7mm}
	\renewcommand\arraystretch{2}
	\caption{Interpretation of $D_2$, $D_1$, $B$, $C$ in likelihood function}
	\begin{tabular}{ccccl}
		\hline
		& $D_2$ & $D_1$ & $B$     & \multicolumn{1}{c}{C} \\ \hline
		Case1 & $m$  & $m$  & $0$     & $ \prod\limits_{i=1}^{m} \left( n-i+1-\sum\limits_{s=1}^{i-1}R_s\right) $ \\
		Case2 & $m$  & $k_1$ & $n-m-\sum\limits_{i=1}^{k_1}R_i$ & $\prod\limits_{i=1}^{m} \left( n-i+1-\sum\limits_{s=1}^{k_1}R_s\right)  \quad$ \\
		Case3 & $k_2$ & $k_1$ & $n-k_2-\sum\limits_{i=1}^{k_1}R_i$ & $\prod\limits_{i=1}^{m} \left( n-i+1-\sum\limits_{s=1}^{k_1}R_s\right)  \quad$ \\ \hline
	\end{tabular}

\label{xishu}
\end{table}

Then, from \eqref{chenfF}, \eqref{L_IAT- II PC}, The likelihood function of $\alpha$ and $\beta$ is given by
	\begin{equation}\label{L}
	\begin{split}
		L(\alpha ,\beta )=&C\left [ \prod_{i=1}^{D_2} \alpha \beta x_i^{\beta -1}\exp \left( \alpha \left( 1-e^{x_i^\beta }\right) +x_i^\beta \right)  \right ] \\
		&\times \left [ \prod_{i=1}^{D_1} \left[ \exp\left(\alpha \left( 1-e^{x_i^\beta }\right) \right) \right] ^{R_i} \right ] \left [  \exp\left( \alpha \left( 1-e^{x_B^\beta }\right) \right)  \right ]^B.
	\end{split}
\end{equation}

Then, we propose several statistical inference methods.

\subsection{Maximum likelihood estimation}

 We can write the natural logarithm of the likelihood function as follows
\begin{equation}\label{lnl}
	\begin{split}
	\ln L(\alpha ,\beta )=&\ln C+D_2\ln\alpha +D_2\ln \beta +\sum_{i=1}^{D_2} (\beta -1)\ln x_i\\&
	+\alpha (1-e^{x_i^\beta })+x_i^\beta+	\sum_{i=1}^{D_1}R_i \alpha (1-e^{x_i^\beta })+B \alpha (1-e^{x_B^\beta }).
\end{split}
\end{equation}

The Maximum likelihood estimations of $\alpha$ and $\beta$ can be obtained by solving the following two equations 
	\begin{align}
\frac{\partial \ln L(\alpha ,\beta )}{\partial \alpha }=&\frac{D_2}{\alpha } +\sum_{i=1}^{D_2}(1-e^{x_i^\beta }) +\sum_{i=1}^{D_1}R_i(1-e^{x_i^\beta })+B(1-e^{x_B^\beta })=0, \label{lalpha} \\
\frac{\partial \ln L(\alpha ,\beta )}{\partial \beta }=&\frac{D_2}{\beta }+\sum_{i=1}^{D_2} \ln x_i+x_i^\beta \ln x_i \notag \\ 
 &-\alpha  \left (\sum_{i=1}^{D_2} e^{x_i^\beta } x_i^\beta  \ln x_i+\sum_{i=1}^{D_1}R_i e^{x_i^\beta } x_i^\beta  \ln x_i+B e^{x_B^\beta }  x_B^\beta  \ln x_B \right )=0 .
\label{lbeta}
\end{align}
From \eqref{lalpha}, the MLE of $\alpha$ can be obtained as follows
\begin{eqnarray}\label{hat_alpha}
	\hat{\alpha } =\frac{D2}{\nu(x_i,\beta )}  ,
\end{eqnarray}
where $\nu(x_i,\beta)=\sum_{i=1}^{D_2}(e^{x_i^\beta }-1) +\sum_{i=1}^{D_1}R_i(e^{x_i^\beta }-1)+B(e^{x_B^\beta }-1)$.
Now substitute equation \eqref{hat_alpha} into equation \eqref{lbeta}, substituting $\hat{\alpha }$ for $\alpha$ to get the MLE estimate of $\beta$ as follows
\begin{eqnarray}\label{fbeta}
   \frac{D_2}{\beta }+\sum_{i=1}^{D_2} \ln x_i(1+x_i^\beta)-\frac{D2}{\nu(x_i,\beta )} \left ( \sum_{i=1}^{D_2} \phi_i (x_i,\beta )+\sum_{i=1}^{D_1}R_i \phi_i(x_i,\beta )+B \phi_B(x_B,\beta ) \right )=0, 
\end{eqnarray}
where $\phi_i (x_i,\beta )=e^{x_i^\beta } x_i^\beta  \ln x_i$.
Equation \eqref{fbeta} does not have a closed-form solution, so the MLE of $\beta$ can be obtained by numerical methods of the following nonlinear equations
\begin{eqnarray}\label{gbeta}
	g(\beta )=\left ( -\frac{\sum_{i=1}^{D_2} \ln x_i(1+x_i^\beta)}{D_2} +\frac{ \sum_{i=1}^{D_2} \phi_i(x_i,\beta )+\sum_{i=1}^{D_1}R_i \phi_i(x_i,\beta )+B \phi_B(x_B,\beta ) }{\nu(x_i,\beta )} \right ) ^{-1}.
\end{eqnarray}

The MLE of $\beta$  is denoted by $\hat{\beta }$, and the solution for $g(\beta )=\beta$ can be obtained by a simple iterative method
\begin{eqnarray*}
	\beta ^{(k+1)} =g(\beta ^{(k)}),
\end{eqnarray*}
here $\beta ^{(k)}$ represents the $k$-$th$ iteration, and the iterative process stops when the difference between two consecutive solutions is less than a predetermined tolerance limit. When we obtain the MLE of $\beta$, the MLE of  $\alpha$  can be obtained directly from equation \eqref{hat_alpha}.
\subsection{ Approximate Confidence Intervals}
In this subsection, we further discuss the $(1-\gamma) \times 100\%$ approximate confidence intervals for the unknown parameters $\lambda=(\alpha,\beta)$ based on the asymptotic distribution of the MLE $\hat{\lambda}$ is $(\hat{\lambda}-\lambda)\to N_2(0,I^{-1}(\lambda))$ (see[\cite{1982}]). where $I^{-1}(\lambda)$ is defined as the inverse Fisher information matrix of the two unknown parameters. \\

The fisher information matrix $I(\lambda )$ is written as which elements are negatives of expected values of the second partial derivatives of the $\ln L(\alpha,\beta)$. The results are as follows
\begin{eqnarray}\label{fisher_information}
	I(\lambda )=-E \begin{bmatrix}\frac{\partial^2\ln L(\alpha ,\beta )}{\partial^2\alpha}
		& \frac{\partial^2\ln L(\alpha ,\beta )}{\partial\alpha \partial\beta }\\\frac{\partial^2\ln L(\alpha ,\beta )}{\partial\beta\partial\alpha }
		&\frac{\partial^2\ln L(\alpha ,\beta )}{\partial^2\beta}
	\end{bmatrix}.
\end{eqnarray}

In general, it can be shown that the asymptotic variance-covariance matrix of $MLE$ is inverting the fisher information matrix $I^{-1}(\lambda )$. The specific details are as follows
\begin{eqnarray}\label{fisher_ni}
	I^{-1}(\lambda)= \begin{bmatrix}-\frac{\partial^2\ln L(\alpha ,\beta )}{\partial^2\alpha}
		& -\frac{\partial^2\ln L(\alpha ,\beta )}{\partial\alpha \partial\beta }\\-\frac{\partial^2\ln L(\alpha ,\beta )}{\partial\beta\partial\alpha }
		&-\frac{\partial^2\ln L(\alpha ,\beta )}{\partial^2\beta}
	\end{bmatrix}  _{(\alpha ,\beta )=(\hat{\alpha} ,\hat{\beta })}^{-1}
	=\begin{bmatrix}var(\hat{\alpha} )
		& cov(\hat{\alpha} ,\hat{\beta } )\\cov(\hat{\beta }  ,\hat{\alpha}) 
		&var(\hat{\beta }  )
	\end{bmatrix}.
\end{eqnarray}

From the log-likelihood function in \eqref{lnl} we can obtain the second derivatives of $\ln L(\alpha ,\beta )$ as follow
\begin{eqnarray}\label{l2}
	 \frac{\partial^2\ln L(\alpha ,\beta )}{\partial^2\alpha }&=&-\frac{D_2}{\alpha ^2},\nonumber\\
	  \frac{\partial^2\ln L(\alpha ,\beta )}{\partial^2\beta }&=&-\frac{D_2}{\beta ^2}+\sum_{i=1}^{D_2} \ln^2x_i x_i^\beta -
	\alpha \left ( \sum_{i=1}^{D_2}\phi_i \xi _i+\sum_{i=1}^{D_1}  R_i \phi_i \xi _i+B\phi_B \xi _B \right ),\\
	 \frac{\partial^2\ln L(\alpha ,\beta )}{\partial\beta \partial \alpha   }&=&\frac{\partial^2\ln L(\alpha ,\beta )}{\partial \alpha\partial\beta \   }=
	-\left ( \sum_{i=1}^{D_2} \phi_i+\sum_{i=1}^{D_1} R_i \phi_i+B \phi_B \right ) ,\nonumber
\end{eqnarray}
where $\xi_i=\ln x_i (1+x_i^{\beta})$,  $\xi _B=\ln x_B (1+x_B^{\beta})$, $\frac{\partial \phi_i (x_i,\beta )}{\partial \beta } =\exp(x_i^{\beta} ) x_i^{\beta}  \ln x_i \xi _i$.

 Therefore, we can obtain the $(1-\gamma) \times 100\%$ approximate confidence intervals for the parameters $\alpha$ and $\beta$ as follow
\begin{eqnarray}\label{confidence_interval}
\left ( {\Large }\hat{\alpha} + Z_{\gamma /2}\sqrt{var(\alpha)},\ \hat{\alpha} + Z_{\gamma /2}\sqrt{var(\alpha)}\right ),\quad
	 \left ( {\Large } \hat{\beta} +Z_{\gamma /2}\sqrt{var(\beta )},\ \hat{\beta} +Z_{\gamma /2}\sqrt{var(\beta )}\right ),
\end{eqnarray}
where $Z_{\gamma /2}$ is the upper $(\gamma /2)^{th}$ percentile point of a standard normal distribution.

\section{Bayesian estimation}
This section is devoted to obtain the Bayes estimators of the parameters $\alpha$ and $\beta$ of Chen distribution based on improved adaptive asymptotic type-II censored data. The Bayesian estimates are obtained using symmetric as well as asymmetric loss function such as squared error loss function(SEL), LINEX loss function(LL) and entropy loss(EL) function. We assume the parameters $\alpha$ and $\beta$ independent and have gamma prior distributions with the following prior distribution. $G(a,b)$ and $G(c,d)$ distribution with PDFs respectively as
\begin{eqnarray*}
	&\pi _1(\alpha ;a,b)=\frac{b^a}{\Gamma (a)}\alpha ^{a-1} e^{-b\alpha },\quad \alpha >0,a,b>0,\\
	&\pi _2(\beta ;c,d)=\frac{d^c}{\Gamma (c)}\beta ^{c-1} e^{-d\beta },\quad \beta >0,c,d>0.	
\end{eqnarray*}

Therefore, the joint prior distribution of  $\alpha$ and $\beta$ is given by
\begin{eqnarray}\label{joint_prior}
	\pi(\alpha ,\beta )\propto \alpha ^{a-1}\beta ^{c-1}e^{-(b\alpha +d\beta )}\quad \alpha ,\beta >0,a,b,c,d>0.
\end{eqnarray}

Subsequently, the joint posterior distribution of  $\alpha$ and $\beta$  becomes
\begin{eqnarray*}\label{posterior}
	\pi (\alpha ,\beta \mid X)=\frac{L(\alpha ,\beta X)\pi (\alpha ,\beta )}{\int_{0}^{\infty } \int_{0}^{\infty }L(\alpha ,\beta X)\pi (\alpha ,\beta )\rm d\alpha \rm d\beta } .
\end{eqnarray*}

Given $\alpha$, $\beta$, directly calculate the joint posterior distribution $\pi (\alpha ,\beta \mid X)$ as follow
\begin{eqnarray}\label{joint_posterior}
	\begin{aligned}
	\pi (\alpha ,\beta \mid X)\propto &\alpha^{D_2+a-1} \exp \left[ -\alpha \left( b-\sum_{i=1}^{D_2}\left( 1-e^{x_i^\beta }\right) -
	\sum_{i=1}^{D_1}R_i\left( 1-e^{x_i^\beta }\right) \right. \right. \\ 
&\left. \left. -B \left( 1-e^{x_B^\beta }\right)  \right)  \right] 
  \times \beta ^{D_2+c-1} \exp\left[ -\beta\left(  d-\sum_{i=1}^{D_2}  logx_i\right)  \right]    \exp\left( \sum_{i=1}^{D_2}x_i^\beta  \right) .	
 \end{aligned}
\end{eqnarray}

 Under the square error loss function (SEL), the Bayes estimate for any parameter $\mu$ is given by
 \begin{eqnarray}\label{SEL}
 	\hat{d}_{SEL}=E(\eta \mid \mu  ) .
 \end{eqnarray}

 For LINEX loss function(LL), the Bayes estimate for any parameter $u$ is given by 
  \begin{eqnarray}\label{LL}
  	\hat{d}_{LL}=-\frac{1}{g}\ln E_\eta (e^{-g\eta} \mid \mu  ) .
  \end{eqnarray}

  For entropy loss function(EL), the Bayesian estimates obtained by minimizing the risk function are
\begin{eqnarray}\label{EL}
	\hat{d}_{EL}=\left [ E_\eta ({\eta^{-q}} \mid \mu  )  \right ] ^{-\frac{1}{q} }.
\end{eqnarray}

Where $g$ and $q$ are known constants, $\eta$ represents any one of the unknown parameters. Obviously, \eqref{SEL} \eqref{LL} \eqref{EL} cannot get an explicit solution. Therefore, the Markov chain Carlo method(MCMC) is used to generate samples from the posterior density function and in turn to compute the Bayesian estimation of the unknown parameters. The following two technique are introduced to calculate bayes estimation.
\subsection{Metropolis-Hastings algorithm with Gibbs sampling}
The Metropolis-Hastings (M-H) algorithm as a general  Markov chain Monte Carlo (MCMC) method be introduced by \cite{Metropolis1953} and then \cite{Hasting1970} extended the M-H algorithm. Gibbs sampling method is a special case of the MCMC method. It can be used to generated a sequence of sampling from the full conditional probability distributions of random variables. Gibbs sampling requires that the joint posterior distribution of each parameter be decomposed into full conditional distribution and then sampling from them. We can apply the Gibbs sampling procedure to generate a sample from \eqref{joint_posterior}, and then compute the Bayesian estimation. \\

Therefore, the posterior conditional density function of $\alpha$ and $\beta$ can be obtained as
\begin{small}
\begin{eqnarray}
	&& \pi (\alpha \mid \beta ,x)\propto \alpha^{D_2+a-1} \exp\left[  -\alpha \left[  b-\sum_{i=1}^{D_2}( 1-e^{x_i^\beta }) -
		\sum_{i=1}^{D_1}R_i( 1-e^{x_i^\beta })  -B ( 1-e^{x_B^\beta }) \right]    \right]   \label{posterior_alpha},\\
	&& \ \pi (\beta \mid \alpha ,x) \propto \beta ^{D_2+c-1} 
		\exp\left[ -\beta \left( d-\sum_{i=1}^{D_2}  \ln x_i\right) \right]   \exp \left( \sum_{i=1}^{D_2}x_i^\beta  \right) \label{posterior_beta},
\end{eqnarray}
\end{small}
here, we use the M-H algorithm with Gibbs sampling to generate samples from \eqref{posterior_alpha} and \eqref{posterior_beta}. The M-H algorithm is performed for Bayesian estimates.\\

\begin{algorithm}[H]
	\caption{ Metropolis-Hastings algorithm} 
	\label{alg1} 
	\begin{algorithmic}[1] 
		\REQUIRE Bayesian estimators of the parameters $\alpha $ and $\beta$ under different loss function with Metropolis-Hastings technique. 
		\ENSURE  \quad \\                                                     
		Step1 : The initial value of the given parameters is $\eta ^{(0)}=(\alpha^{(0)},\beta^{(0)})$.\\
		Step2 : Use Metropolis-Hasting algorithm to generate $\alpha ^{(0)}$ from\eqref{posterior_alpha} and $\beta ^{(0)}$ from \eqref{posterior_beta} as the Gibbs sampling iteration values. $q(\beta)=N(\hat{\beta})$. Take the normal distribution as the  proposal distribution. Repeat $N$ times, then we get $\eta ^{(h)}=(\alpha^{(h)},\beta^{(h)}),\quad h=1,2,...,N$.\\
		Step3 : Take the first $M$ results as a burn-in phase. Under different loss function, Bayesian estimations are 
		\begin{eqnarray*}\label{hat_loss}
			\hat{d}_{SEL}&=&\frac{1}{N-M}\sum_{h=M+1}^{N}\eta^{(h)},\\
			\hat{d}_{LL}&=&-\frac{1}{g}\ln\left(  \frac{1}{N-M}\sum_{h=M+1}^{N}e^{g\eta^{{\scriptsize (} h{\scriptsize )} }}\right) ,\\
			\hat{d}_{EL}&= & \left[\frac{1}{N-M}  \sum_{h=M+1}^{N} \left( \eta^{{\scriptsize (} h{\scriptsize )} } \right) ^{-q} \right ] ^{-\frac{1}{q} }.	
		\end{eqnarray*}
		Step4 : Get the results in step2. Sort in ascending order as $\left( \eta_{(1)},\eta_{(2)},...,\eta_{(N)}\right) $.\\
	\end{algorithmic} 
\end{algorithm}

\subsection{Importance sampling technique}
The importance sampling method can be used to compute the approximate results of \eqref{SEL}, \eqref{LL} and \eqref{EL}, so as to obtain the Bayesian estimations.
From \eqref{joint_posterior}, we have 
\begin{eqnarray}\label{pi}
	\pi (\alpha ,\beta \mid X)\propto  g(\alpha \mid \beta ,x)  g(\beta \mid x) \omega (\alpha ,\beta ),
\end{eqnarray}
where
\begin{eqnarray}\label{pif}
 	g(\alpha \mid \beta ,x)&\propto &\alpha^{D_2+a-1} \exp\left[ -\alpha \left[   b+
		\sum_{i=1}^{D_1}R_i\left( e^{x_i^\beta }-1\right)  +B \left( e^{x_B^\beta }-1\right) \right]   \right] ,\\
	g(\beta \mid x)&\propto &\beta ^{D_2+c-1} \exp\left[ -\beta \left( d-\sum_{i=1}^{D_2}  logx_i\right) \right] ,\\
	\omega (\alpha ,\beta)&\propto & \frac{\exp\left[ \sum_{i=1}^{D_2} \alpha \left( 1-e^{x_i^\beta }\right) +\sum_{i=1}^{D_2}x_i^\beta  \right] }{ \left [ b+\nu\left( x_i,\beta\right)  \right ]  ^{D_2+a}},    \  \  \   \  \ 
\end{eqnarray}
where  $\nu(x_i,\beta)=\sum_{i=1}^{D_2}(e^{x_i^\beta }-1) +\sum_{i=1}^{D_1}R_i(e^{x_i^\beta }-1)+B(e^{x_B^\beta }-1)$.
Note that the distribution of $g(\alpha \mid \beta, x)$ follows Gamma distribution with parameters$(D_2+a-1)$ and $(b-\sum_{i=1}^{D_2}(1-e^{x_i^\beta })-
\sum_{i=1}^{D_1}R_i(1-e^{x_i^\beta }) -B (1-e^{x_B^\beta }))$. The distribution of  $g(\beta \mid x)$ follow Gamma distribution with parameters $(D_2+c-1)$ and $(d-\sum_{i=1}^{D_2}  \log x_i)$. Therefore, one can easily generate samples from the distribution of $\alpha$ and $\beta$, respectively.
Then by the importance sampling method, the Bayesian estimation of $\alpha$ and $\beta$ under the squared error loss function are obtained as
\begin{eqnarray}\label{sel_alpha}
	\hat{\alpha }_{SEL}=\frac{\sum_{i=1}^{N_2}\alpha ^{(i)} \omega ^{(i)}}{\sum_{i=1}^{N_2}\omega ^{(i)}} ,
\end{eqnarray}
and
\begin{eqnarray}\label{sel_beta}
	\hat{\beta}_{SEL}=\frac{\sum_{i=1}^{N_2}\beta ^{(i)} \omega ^{(i)}}{\sum_{i=1}^{N_2}\omega ^{(i)}} .
\end{eqnarray}

Similarly, the Bayesian estimation of  $\alpha$ and $\beta$ under the LINEX loss function are obtained as
\begin{eqnarray}\label{ll_alpha}
	\hat{\alpha}_{LL}=-\frac{1}{g} \ln \left [ \frac{\sum_{i=1}^{N_2}e^{-g\alpha ^{(i)}} \omega ^{(i)}}{\sum_{i=1}^{N_2}\omega ^{(i)}}  \right ] ,
\end{eqnarray}
and
\begin{eqnarray}\label{ll_beta}
	\hat{\beta}_{LL}=-\frac{1}{g} \ln \left [ \frac{\sum_{i=1}^{N_2}e^{-g\beta ^{(i)}} \omega ^{(i)}}{\sum_{i=1}^{N_2}\omega ^{(i)}}  \right ] .
\end{eqnarray}

Again the Bayesian estimation of $\alpha$ and $\beta$ under the entropy loss function are obtained as
\begin{eqnarray}\label{el_alpha}
	\hat{\alpha}_{EL}=\left [ \frac{\sum_{i=1}^{N_2}(\alpha ^{{\scriptsize (} i{\scriptsize )} })^{-q} \omega ^{(i)}}{\sum_{i=1}^{N_2}\omega ^{(i)}}  \right ]^{-\frac{1}{q} } ,
\end{eqnarray}
and 
\begin{eqnarray}\label{el_beta}
	\hat{\beta}_{EL}=\left [ \frac{\sum_{i=1}^{N_2}(\beta ^{{\scriptsize (} i{\scriptsize )} })^{-q} \omega ^{(i)}}{\sum_{i=1}^{N_2}\omega ^{(i)}}  \right ]^{-\frac{1}{q} } .
\end{eqnarray}

The procedure of importance sampling is given as follows\\
\begin{algorithm} 
	\caption{Importance sampling algorithm} 
	\label{alg1} 
	\begin{algorithmic}[2] 
		\REQUIRE  Bayesian estimators of the parameters $\alpha $ and $\beta$ under different loss function with importance sampling technique.
		\ENSURE \quad \\   
	Step1: Generate $\beta^{(1)}$ from $g(\beta \mid x)$.\\
	Step2: Generate $\alpha^{(1)}$ from $g(\alpha \mid \beta ,x)$.\\
	Step3: Calculate $\omega^{(1)} (\beta^{(1)} )$ and  $\theta^{(1)}=\theta(\alpha^{(1)},\beta^{(1)}) $.\\
	Step4: Respect the above process $N_2$ times and obtain $\left(  \omega^{(1)},\omega^{(2)},...,\omega^{(N_2)}\right) $ and $\left( \theta^{(1)},\theta^{(2)},...,\theta^{(N_2)}\right) $.\\
	Step5: Calculate the Bayesian estimators of the parameters $\alpha $ and $\beta$ under different loss function. 
	\end{algorithmic} 
\end{algorithm}

\section{Simulation experiments}
In this section, the simulation study is carried out Monte Carlo simulation to evaluate the performance of the proposed methods. The Biases and mean square errors(MSE) of the MLE and Bayesian estimates for $\alpha$ and $\beta$ are evaluated under the IAT-II PCS and and in $R$ program. The Bayesian estimates are computed by using Importance sampling and M-H technique. Under different combination of $(n,m)$ and different censoring schemes $(T_1,T_2,R_1,R_2,\cdots,R_m)$. the point and interval estimation results are assessed on the basis of mean Bias and MSE respectively.

All the estimates are to compute arbitrary the unknown parameter values $\alpha=0.2$ and $\beta=0.5)$. Accordingly hyperparameters in gamma prior are assigned as $a=2$, $b=2,$ and $c=2$, $d=2$. Here we consider four different censoring schemes, namely:\\
\begin{eqnarray*}
 \textup{Scheme1}&:&R_m=n-m,R_i=0  \ \ \textup{for} \  i\ne m.\\
 \textup{Scheme2}&:&R_1=n-m,R_i=0 \ \ \textup{for}  \ i\ne 1.\\
 \textup{Scheme3}&:&R_{\frac{m+1}{2}}=n-m,  R_i=0 \ \ \textup{for}\   i  \ne\frac{m+1}{2};\ \textup{if}\  m\   \textup{is}\  \textup{odd},\ \textup{and}  \\
        & &  R_{\frac{m}{2}}=n-m,  R_i=0 \ \ \textup{for} \   i\ne\frac{m}{2},\  \textup{if} \ m \  \textup{is}\  \textup{even}.\\
 \textup{Scheme4}&:&R_i=\frac{n-m}{m},\quad i=1,2,\cdots,m. \\
 \end{eqnarray*}

  In our study, two expected time on test $(T_1,T_2)$ are (0.4,4) and (1,7). We consider three different values of $(n,m)$, namely (15,5), (20,10) and (30,15). For evaluating Bayes estimators under LINEX loss function, we take $g=1$ and for entropy loss function, we take $q=1$. Based on the simulation study we have the following conclusion. 
  
  From the simulate results in Table \ref{mle1} and \ref{mle2}, one can observe following conclusions for MLE. From these Tables the following conclusion are made:
  
  1. The estimation of MLE is evaluated by Bias and MSE, among which MSE is better than Bias.
  
  2. In most cases, for different censoring schemes, Scheme 4 is uniform censoring scheme and others are non-uniform schemes. It can be seen from the table that the estimation effect of uniform censoring scheme is better than that of other non-uniform censoring schemes.
  
  3. Compared with Table \ref{mle1} and Table \ref{mle2}, in most cases, when the time threshold becomes larger, the Bias and MSE of related MLE becomes larger, the estimation effect of $(T_1,T_2)=(1,7)$ is poor. Therefore, it is essential to set a reasonable time threshold in the experiment.

 \begin{table}[!ht]
 	\centering
 \caption{{\scriptsize Bias and MSE of the MLE $\hat{\alpha}$ and $\hat{\beta}$ for different choices of $n,m$ and $T_1=0.4,T_2=4$.}}
 \begin{tabular}{ccllrrlrr}
 	\hline
              & \multicolumn{2}{c}{}               & \multicolumn{1}{c}{}     & \multicolumn{2}{c}{$\hat{\alpha}$} &  & \multicolumn{2}{c}{$\hat{\beta}$} \\ \cline{5-6} \cline{8-9} 
 	n,m       & scheme & \multicolumn{1}{c}{$\alpha$} & \multicolumn{1}{c}{$\beta$}& \multicolumn{1}{c}{Bias}       & \multicolumn{1}{c}{MSE}        &  &\multicolumn{1}{c}{Bias}         &  \multicolumn{1}{c}{MSE}        \\ \hline
 	& I      & 0.2                       & 0.5                      & -0.18276     & 0.03340    &  & 0.11049     & 0.01221    \\
 	& II     & 0.2                       & 0.5                      & -0.06933     & 0.00481    &  & -0.24902    & 0.06201    \\
 	15,5  & III    & 0.2                       & 0.5                      & -0.19048     & 0.03628    &  & 0.17394     & 0.03026    \\
 	& IV     & 0.2                       & 0.5                      & -0.03561     & 0.00127    &  & -0.01553    & 0.00024    \\
 	& I      & 0.2                       & 0.5                      & -0.08918     & 0.00795    &  & 0.25224     & 0.06362    \\
 	& II     & 0.2                       & 0.5                      & 0.09200      & 0.00846    &  & -0.23686    & 0.05610    \\
 	20,10 & III    & 0.2                       & 0.5                      & 0.06425      & 0.00413    &  & -0.14137    & 0.01998    \\
 	& IV     & 0.2                       & 0.5                      & -0.05793     & 0.00336    &  & 0.11998     & 0.01440    \\
 	& I      & 0.2                       & 0.5                      & -0.17254     & 0.02977    &  & 0.41297     & 0.17055    \\
 	& II     & 0.2                       & 0.5                      & -0.11648     & 0.01357    &  & -0.35441    & 0.12561    \\
 	30,15 & III    & 0.2                       & 0.5                      & -0.10793     & 0.01165    &  & -0.08880    & 0.00789    \\
 	& IV     & 0.2                       & 0.5                      & 0.02276      & 0.00052    &  & -0.21521    & 0.04632    \\ \hline
 	\end{tabular}
 	\label{mle1}
 \end{table}
 
\begin{table}[!ht]
	\centering
\caption{{\scriptsize Bias and MSE of the MLE $\hat{\alpha}$ and $\hat{\beta}$ for different choices of $n,m$ and $T_1=1,T_2=7$.}}
\begin{tabular}{ccllrrlrr}
	\hline
	         & \multicolumn{2}{c}{}               & \multicolumn{1}{c}{}     & \multicolumn{2}{c}{$\hat{\alpha}$} &  & \multicolumn{2}{c}{$\hat{\beta}$} \\ \cline{5-6} \cline{8-9} 
	n,m       & scheme & \multicolumn{1}{c}{$\alpha$} & \multicolumn{1}{c}{$\beta$} & \multicolumn{1}{c}{Bias}         & \multicolumn{1}{c}{MSE}        &  & \multicolumn{1}{c}{Bias}        & \multicolumn{1}{c}{MSE}        \\ \hline
	                        & I                          & 0.2   & 0.5  & -0.09838     & 0.00968    &  & -0.03055    & 0.00093    \\ 
	& II                         & 0.2   & 0.5  & -0.11620     & 0.01350    &  & -0.34789    & 0.12103    \\
	15,5                & III                        & 0.2   & 0.5  & -0.08268     & 0.00684    &  & -0.34844    & 0.12141    \\
	& IV                         & 0.2   & 0.5  & -0.11235     & 0.01262    &  & 0.21203     & 0.04496    \\
	& I                          & 0.2   & 0.5  & -0.11840     & 0.01402    &  & 0.18201     & 0.03313    \\
	& II                         & 0.2   & 0.5  & 0.12462      & 0.01553    &  & -0.22917    & 0.05252    \\
	20,10               & III                        & 0.2   & 0.5  & -0.11064     & 0.01224    &  & 0.24352     & 0.05930    \\
	& IV                         & 0.2   & 0.5  & -0.12914     & 0.01668    &  & 0.31077     & 0.09658    \\
	& I                          & 0.2   & 0.5  & -0.11563     & 0.01337    &  & 0.19582     & 0.03834    \\
	& II                         & 0.2   & 0.5  & -0.11043     & 0.01219    &  & 0.23403     & 0.05477    \\
	30,15               & III                        & 0.2   & 0.5  & 0.07763      & 0.00603    &  & -0.20020    & 0.04008    \\
	& IV                         & 0.2   & 0.5  & -0.08184     & 0.00670    &  & 0.16316     & 0.02662    \\ \hline
	\end{tabular}
	\label{mle2}
\end{table}

 Table \ref{aci1} and \ref{aci2} present the coverage probability and average length of two-side $95.5\%$ confidence intervals of parameters constructed by maximum likelihood method (using fisher observation information matrix). From these Tables the following conclusion are made:
 
 1. The truth value of $\alpha=0.2$ and $\beta=0.5$ is in the middle of the each interval and can be well covered.
 
 2. The asymptotic confidence interval of $\beta$ obtained by fisher observation information matrix, in most case, the average length of the interval is getting longer when the time threshold increase. Meanwhile, different censoring schemes $R$ had no significant effect on the $\alpha$ and $\beta$ results of interval estimation.
 
 3. The average length of the interval about $\alpha$ is significantly different under the two time thresholds $(0.4,4)$ and $(1,7)$. 
 
 4. The asymptotic confidence interval of $\alpha$, in most case, are more ideal than that of $\beta$.

\begin{table}[H]
	\centering
	\caption{{\scriptsize Interval estimations,and Average Length(AL) for different choices of $n,m$ and $T_1=0.4,T_2=4$.}}
	\begin{tabular}{cccccc}
		\hline
		n,m & scheme & $\alpha$ & $\beta$ & $\alpha-AL$ & $\beta-AL$ \\ \hline
		~ & I & (0.00129,0.58393) & (0.48694,1.71175) & 0.58264  & 1.22482  \\ 
		~ & II & (0.02411,0.28934) & (0.04585,0.63700) & 0.26523  & 0.59115  \\ 
		15,5 & III & (0.02096,0.40454) & (0.11858,1.42586) & 0.38358  & 1.30727  \\ 
		~ & IV & (0.06947,0.71424) & (0.09405,0.86234) & 0.64476  & 0.76829  \\ 
		~ & I & (0.01640,0.70959) & (0.17385,0.90605) & 0.57258  & 0.73220  \\ 
		~ & II & (0.13701,0.53979) & (0.01334,0.04607) & 0.38361  & 0.03273  \\ 
		20,10 & III & (0.13223,0.75290) & (0.14972,0.85257) & 0.62068  & 0.70285  \\ 
		~ & IV & (0.10312,0.54511) & (0.14693,0.65170) & 0.44199  & 0.50477  \\ 
		~ & I & (0.03949,0.31102) & (0.45257,1.11968) & 0.27154  & 0.66710  \\ 
		~ & II & (0.15002,0.56836) & (0.11233,0.51236) & 0.41834  & 0.40002  \\ 
		30,15 & III & (0.12563,0.55586) & (0.16401,0.78419) & 0.43023  & 0.62018  \\ 
		~ & IV & (0.06050,0.33362) & (0.36250,1.06036) & 0.27311  & 0.69786 \\ \hline
	\end{tabular}
	\label{aci1}
\end{table}

\begin{table}[H]
	\centering
	\caption{{\scriptsize Interval estimations,and Average Length(AL) for different choices of $n,m$ and $T_1=1,T_2=7$.}}
	\begin{tabular}{cccccc}
		\hline
			n,m & scheme & $\alpha$ & $\beta$ & $\alpha-AL$ & $\beta-AL$ \\ \hline 
		~ & I & (0.02227,0.46365) &(0.17421,1.26507) & 0.44138  &1.09087\\	
		~ & II & (0.02410,0.29141) & (0.04789,0.66570) & 0.26731  & 0.61781  \\ 
		15,5 & III & (0.02988,0.46066) & (0.04098,0.80587) & 0.43079  & 0.76489  \\ 
		~ & IV & (0.02385,0.32205) & (0.31173,2.04524) & 0.29820  & 1.73351  \\ 
		~ & I & (0.02728,0.73341) & (0.09402,0.61997) & 0.58973  & 0.52595  \\ 
		~ & II & (0.13713,0.76841) & (0.08804,0.65271) & 0.63127  & 0.56467  \\ 
		20,10 & III & (0.02880,0.27724) & (0.46691,1.07021) & 0.24844  & 0.60330  \\ 
		~ & IV & (0.01690,0.29710) & (0.47157,1.23919) & 0.28020  & 0.76762  \\ 
		~ & I & (0.03392,0.20986) & (0.48876,0.96220) & 0.17594  & 0.47344  \\ 
		~ & II & (0.04174,0.33184) & (0.38175,1.41139) & 0.29011  & 1.02963  \\ 
		30,15 & III & (0.14169,0.54398) & (0.15428,0.58258) & 0.40229  & 0.40229  \\ 
		~ & IV & (0.05112,0.27315) & (0.41602,1.05710) & 0.22203  & 0.64107 \\ \hline
	\end{tabular}
	\label{aci2}
\end{table}

Table \ref{bayesa1}, \ref{bayesb1}, \ref{bayesa2}, \ref{bayesb2} present the Bias and MSE of Bayesian estimation of $\alpha$ and $\beta$ by importance sampling technique under squared error loss function, LINEX loss function and entropy loss function. Table \ref{bayesa3}, \ref{bayesb3}, \ref{bayesa4}, \ref{bayesb4} present the Bias and MSE of Bayesian estimation by Metropolis-Hasting technique of $\alpha$ and $\beta$. From these Tables the following conclusion are made:

1. For the tables, we can also see that Bias and MSE of Bayes estimators of $\alpha$ and $\beta$ are smaller than Bias and MSE of corresponding MLE. 

2. From Table \ref{bayesa1}, \ref{bayesa2} and Table \ref{bayesb1}, \ref{bayesb2}, in most cases, one can see that Bayesian estimators of $\alpha$ and $\beta$ by importance sampling technique under squared error loss function possess minimum Bias and MSE. When the time threshold increase, under entropy loss function also possess minimum Bias and MSE. From Table \ref{bayesa3} and Table \ref{bayesa4}, one can see Bayesian estimators of $\alpha$ by M-H technique under squared error loss function possess minimum Bias and MSE. From Table \ref{bayesb3} Bayesian estimators of $\beta$ by M-H technique under LINEX loss function possess minimum Bias and MSE, when the time threshold is $(T_1,T_2)=(0.4,4)$ and Table \ref{bayesb4}, Bayesian estimators of $\beta$ by M-H technique under LINEX loss function and entropy loss function possess minimum Bias and MSE, when the time threshold is $(T_1,T_2)=(1,7)$.

3. The Bias and MSE of Bayesian estimation using M-H technique, in most case, are smaller than that using importance sampling.

\begin{table}[H]
	\centering
	\caption{{\scriptsize  Bias and MSE of the Bayesian estimates by importance sampling technique of $\hat{\alpha}$ under squared error loss function $\hat{\alpha}_{SEL}$, LINEX loss function $\hat{\alpha}_{LL}$, and entropy loss function $\hat{\alpha}_{EL}$, for different choices of $n, m,\alpha,\beta$ and $T_1 = 0.4,T_2=4$.}}
\begin{tabular}{ccccrrrrrrrr}
	\hline
	&     &       &      & \multicolumn{2}{c}{$\hat{\alpha}_{SEL}$} &  & \multicolumn{2}{c}{$\hat{\alpha}_{LL}$} &  & \multicolumn{2}{c}{$\hat{\alpha}_{EL}$} \\ \cline{5-6} \cline{8-9} \cline{11-12} 
	n,m   & sc  & alpha & beta & \multicolumn{1}{c}{Bias}         & \multicolumn{1}{c}{MSE}         &  & \multicolumn{1}{c}{Bias}        & \multicolumn{1}{c}{MSE}        &  & \multicolumn{1}{c}{Bias}          & \multicolumn{1}{c}{MSE}         \\ \hline
	& I   & 0.2   & 0.5  & -0.05442     & 0.00296     &  & -0.05454     & 0.00297    &  & -0.05571      & 0.00310     \\
	& II  & 0.2   & 0.5  & 0.00409      & 0.00002     &  & 0.00236      & 0.00001    &  & -0.01083      & 0.00012     \\
	15,5  & III & 0.2   & 0.5  & -0.03628     & 0.00132     &  & -0.03630     & 0.00132    &  & -0.03644      & 0.00133     \\
	& IV  & 0.2   & 0.5  & -0.08628     & 0.00744     &  & -0.08746     & 0.00765    &  & -0.10018      & 0.01004     \\
	& I   & 0.2   & 0.5  & -0.05352     & 0.00286     &  & -0.05365     & 0.00288    &  & -0.05504      & 0.00303     \\
	& II  & 0.2   & 0.5  & -0.07435     & 0.00553     &  & -0.07471     & 0.00558    &  & -0.07949      & 0.00632     \\
	20,10 & III & 0.2   & 0.5  & 0.01078      & 0.00012     &  & 0.01070      & 0.00011    &  & 0.01009       & 0.00010     \\
	& IV  & 0.2   & 0.5  & 0.01141      & 0.00013     &  & 0.00952      & 0.00009    &  & -0.00244      & 0.00001     \\
	& I   & 0.2   & 0.5  & -0.10540     & 0.01111     &  & -0.10750     & 0.01156    &  & -0.13787      & 0.01901     \\
	& II  & 0.2   & 0.5  & 0.05013      & 0.00251     &  & 0.04831      & 0.00233    &  & 0.03676       & 0.00135     \\
	30,15 & III & 0.2   & 0.5  & -0.02612     & 0.00068     &  & -0.02675     & 0.00072    &  & -0.03156      & 0.00100     \\
	& IV  & 0.2   & 0.5  & 0.01169      & 0.00014     &  & 0.01142      & 0.00013    &  & 0.00932       & 0.00009     \\ \hline
	\end{tabular}
	\label{bayesa1}
\end{table}
\begin{table}[H]
	\centering
	\caption{{\scriptsize Bias and MSE of the Bayesian estimates by importance sampling technique of $\hat{\beta}$ under squared error loss function $\hat{\beta}_{SEL}$ LINEX loss function $\hat{\beta}_{LL}$ and entropy loss function  $\hat{\beta}_{EL}$ for different choices of $n, m,\alpha,\beta$ and $T_1 = 0.4,T_2=4$.}}
	\begin{tabular}{ccccrrrrrrrr}
			\hline
			& \multicolumn{2}{c}{} &      & \multicolumn{2}{c}{$\hat{\beta}_{SEL}$} &  & \multicolumn{2}{c}{$\hat{\beta}_{LL}$} &  & \multicolumn{2}{c}{$\hat{\beta}_{EL}$} \\ \cline{5-6} \cline{8-9} \cline{11-12} 
			n,m       & sc    & $\alpha$    & $\beta$ & \multicolumn{1}{c}{Bias}         & \multicolumn{1}{c}{MSE}         &  & \multicolumn{1}{c}{Bias}         & \multicolumn{1}{c}{MSE}        &  & \multicolumn{1}{c}{Bias}          & \multicolumn{1}{c}{MSE}         \\ \hline
& I         & 0.2      & 0.5  & 0.06578      & 0.00433     &  & 0.06322      & 0.00400    &  & 0.05649       & 0.00319     \\
& II        & 0.2      & 0.5  & 0.11644      & 0.01356     &  & 0.10269      & 0.01055    &  & 0.04396       & 0.00193     \\
15,5  & III       & 0.2      & 0.5  & 0.09152      & 0.00838     &  & 0.09152      & 0.00838    &  & 0.09152       & 0.00838     \\
& IV        & 0.2      & 0.5  & 0.07352      & 0.00541     &  & 0.06354      & 0.00404    &  & 0.04380       & 0.00192     \\
& I         & 0.2      & 0.5  & -0.04531     & 0.00205     &  & -0.04584     & 0.00210    &  & -0.04785      & 0.00229     \\
& II        & 0.2      & 0.5  & 0.15521      & 0.02409     &  & 0.12921      & 0.01670    &  & 0.01849       & 0.00034     \\
20,10 & III       & 0.2      & 0.5  & -0.00909     & 0.00008     &  & -0.00917     & 0.00008    &  & -0.00943      & 0.00009     \\
& IV        & 0.2      & 0.5  & -0.01205     & 0.00015     &  & -0.01630     & 0.00027    &  & -0.03579      & 0.00128     \\
& I         & 0.2      & 0.5  & 0.21632      & 0.04679     &  & 0.21328      & 0.04549    &  & 0.20346       & 0.04140     \\
& II        & 0.2      & 0.5  & 0.25439      & 0.06471     &  & 0.24590      & 0.06047    &  & 0.21985       & 0.04834     \\
30,15 & III       & 0.2      & 0.5  & -0.07071     & 0.00500     &  & -0.07143     & 0.00510    &  & -0.07462      & 0.00557     \\
& IV        & 0.2      & 0.5  & 0.15959      & 0.02547     &  & 0.15892      & 0.02526    &  & 0.15735       & 0.02476     \\ \hline
	\end{tabular}
	\label{bayesb1}
\end{table}

\begin{table}[H]
	\centering
	\caption{{\scriptsize Bias and MSE of the Bayesian estimates by importance sampling technique of $\hat{\alpha}$ under squared error loss function $\hat{\alpha_{SEL}}$, LINEX loss function $\hat{\alpha_{LL}}$, and entropy loss function $\hat{\alpha_{EL}}$, for different choices of $n, m,\alpha,\beta$ and $T_1 = 1,T_2=7$.}}
	\begin{tabular}{ccccrrrrrrrr}
		\hline
		& \multicolumn{2}{c}{} &      & \multicolumn{2}{c}{$\hat{\alpha_{SEL}}$} &  & \multicolumn{2}{c}{$\hat{\alpha_{LL}}$} &  & \multicolumn{2}{c}{$\hat{\alpha_{EL}}$} \\ \cline{5-6} \cline{8-9} \cline{11-12} 
		n,m       & sc    & $\alpha$    & $\beta$ & \multicolumn{1}{c}{Bias}         & \multicolumn{1}{c}{MSE}         &  & \multicolumn{1}{c}{Bias}         & \multicolumn{1}{c}{MSE}        &  & \multicolumn{1}{c}{Bias}          & \multicolumn{1}{c}{MSE}         \\ \hline
		& I         & 0.2      & 0.5  & 0.02854      & 0.00081     &  & 0.02854      & 0.00081    &  & 0.02845       & 0.00081     \\
		& II        & 0.2      & 0.5  & -0.05099     & 0.00260     &  & -0.05099     & 0.00260    &  & -0.05576      & 0.00311     \\
		15,5  & III       & 0.2      & 0.5  & 0.06638      & 0.00441     &  & 0.06638      & 0.00441    &  & 0.06353       & 0.00404     \\
		& IV        & 0.2      & 0.5  & 0.04207      & 0.00177     &  & 0.04207      & 0.00177    &  & 0.04180       & 0.00175     \\
		& I         & 0.2      & 0.5  & -0.00651     & 0.00004     &  & -0.00651     & 0.00004    &  & -0.01564      & 0.00024     \\
		& II        & 0.2      & 0.5  & 0.05089      & 0.00259     &  & 0.05089      & 0.00259    &  & 0.04696       & 0.00220     \\
		20,10 & III       & 0.2      & 0.5  & -0.09754     & 0.00951     &  & -0.09754     & 0.00951    &  & -0.09998      & 0.01000     \\
		& IV        & 0.2      & 0.5  & 0.09069      & 0.00823     &  & 0.09069      & 0.00823    &  & 0.08205       & 0.00673     \\
		& I         & 0.2      & 0.5  & -0.05424     & 0.00294     &  & -0.05424     & 0.00294    &  & -0.05457      & 0.00298     \\
		& II        & 0.2      & 0.5  & -0.06160     & 0.00379     &  & -0.06160     & 0.00379    &  & -0.06420      & 0.00412     \\
		30,15 & III       & 0.2      & 0.5  & -0.04417     & 0.00195     &  & -0.04417     & 0.00195    &  & -0.04914      & 0.00241     \\
		& IV        & 0.2      & 0.5  & 0.00413      & 0.00002     &  & 0.00413      & 0.00002    &  & 0.00167       & 0.00000     \\ \hline
	\end{tabular}
	\label{bayesa2}
\end{table}

\begin{table}[H]
	\centering
	\caption{{\scriptsize Bias and MSE of the Bayesian estimates by importance sampling technique of $\hat{\beta}$ under squared error loss function $\hat{\beta}_{SEL}$ LINEX loss function $\hat{\beta}_{LL}$ and entropy loss function  $\hat{\beta}_{EL}$ for different choices of $n, m,\alpha,\beta$ and $T_1 = 1,T_2=7$.}}
\begin{tabular}{ccccrrrrrrrr}
	\hline
	& \multicolumn{2}{c}{} &      & \multicolumn{2}{c}{$\hat{\beta}_{SEL}$} &  & \multicolumn{2}{c}{$\hat{\beta}_{LL}$} &  & \multicolumn{2}{c}{$\hat{\beta}_{EL}$} \\ \cline{5-6} \cline{8-9} \cline{11-12} 
	n,m       & sc   & $\alpha$    & $\beta$ & \multicolumn{1}{c}{Bias}         & \multicolumn{1}{c}{MSE}         &  & \multicolumn{1}{c}{Bias}         & \multicolumn{1}{c}{MSE}        &  & \multicolumn{1}{c}{Bias}          & \multicolumn{1}{c}{MSE}         \\ \hline
	& I         & 0.2      & 0.5  & -0.01712     & 0.00029     &  & -0.01712     & 0.00029    &  & -0.01714      & 0.00029     \\
	& II        & 0.2      & 0.5  & -0.15450     & 0.02387     &  & -0.15542     & 0.02416    &  & -0.16039      & 0.02572     \\
	15,5  & III       & 0.2      & 0.5  & -0.10080     & 0.01016     &  & -0.30124     & 0.09074    &  & -0.30434      & 0.09262     \\
	& IV        & 0.2      & 0.5  & 0.06374      & 0.00406     &  & -0.13645     & 0.01862    &  & -0.13695      & 0.01876     \\
	& I         & 0.2      & 0.5  & 0.08405      & 0.00706     &  & 0.07961      & 0.00634    &  & 0.04465       & 0.00199     \\
	& II        & 0.2      & 0.5  & 0.02863      & 0.00082     &  & 0.02693      & 0.00073    &  & 0.01339       & 0.00018     \\
	20,10 & III       & 0.2      & 0.5  & 0.03779      & 0.00143     &  & 0.03532      & 0.00125    &  & 0.02682       & 0.00072     \\
	& IV        & 0.2      & 0.5  & -0.11941     & 0.01426     &  & -0.12013     & 0.01443    &  & -0.12269      & 0.01505     \\
	& I         & 0.2      & 0.5  & 0.07355      & 0.00541     &  & 0.07316      & 0.00535    &  & 0.07221       & 0.00521     \\
	& II        & 0.2      & 0.5  & 0.04997      & 0.00250     &  & 0.04939      & 0.00244    &  & 0.04749       & 0.00226     \\
	30,15 & III       & 0.2      & 0.5  & -0.06077     & 0.00369     &  & -0.06093     & 0.00371    &  & -0.06154      & 0.00379     \\
	& IV        & 0.2      & 0.5  & 0.06078      & 0.00369     &  & 0.06014      & 0.00362    &  & 0.05806       & 0.00337     \\ \hline
	\end{tabular}
	\label{bayesb2}
\end{table}

\begin{table}[H]
	\centering
	\caption{{\scriptsize Bias and MSE of the Bayesian estimates by M-H technique of $\hat{\alpha}$ under squared error loss function $\hat{\alpha}_{SEL}$, LINEX loss function $\hat{\alpha}_{LL}$, and entropy loss function $\hat{\alpha}_{EL}$, for different choices of $n, m,\alpha,\beta$ and $T_1 = 0.4,T_2=4$.}}
\begin{tabular}{ccccrrrrrrrr}
	\hline
	& \multicolumn{2}{c}{} &      & \multicolumn{2}{c}{$\hat{\alpha}_{SEL}$} &  & \multicolumn{2}{c}{$\hat{\alpha}_{LL}$} &  & \multicolumn{2}{c}{$\hat{\alpha}_{EL}$} \\ \cline{5-6} \cline{8-9} \cline{11-12} 
	n,m       & sc   & $\alpha$    & $\beta$ & \multicolumn{1}{c}{Bias}         & \multicolumn{1}{c}{MSE}         &  & \multicolumn{1}{c}{Bias}         & \multicolumn{1}{c}{MSE}        &  & \multicolumn{1}{c}{Bias}          & \multicolumn{1}{c}{MSE}         \\ \hline
	& I         & 0.2      & 0.5  & -0.03289     & 0.00108     &  & -0.05811     & 0.00338    &  & -0.16499      & 0.02722     \\
	& II        & 0.2      & 0.5  & 0.02562      & 0.00066     &  & 0.01058      & 0.00011    &  & -0.00766      & 0.00006     \\
	15,5  & III       & 0.2      & 0.5  & -0.03411     & 0.00116     &  & -0.07142     & 0.00510    &  & -0.13038      & 0.01700     \\
	& IV        & 0.2      & 0.5  & 0.05093      & 0.00259     &  & 0.01892      & 0.00036    &  & -0.03241      & 0.00105     \\
	& I         & 0.2      & 0.5  & -0.02719     & 0.00074     &  & -0.07549     & 0.00570    &  & -0.16503      & 0.02723     \\
	& II        & 0.2      & 0.5  & -0.00444     & 0.00002     &  & -0.02316     & 0.00054    &  & -0.03929      & 0.00154     \\
	20,10 & III       & 0.2      & 0.5  & 0.00170      & 0.00000     &  & -0.03327     & 0.00111    &  & -0.07890      & 0.00622     \\
	& IV        & 0.2      & 0.5  & -0.01749     & 0.00031     &  & -0.04216     & 0.00178    &  & -0.14863      & 0.02209     \\
	& I         & 0.2      & 0.5  & 0.08353      & 0.00698     &  & 0.06098      & 0.00372    &  & -0.03633      & 0.00132     \\
	& II        & 0.2      & 0.5  & -0.03041     & 0.00092     &  & -0.05154     & 0.00266    &  & -0.11496      & 0.01322     \\
	30,15 & III       & 0.2      & 0.5  & 0.06115      & 0.00374     &  & 0.04258      & 0.00181    &  & 0.01858       & 0.00035     \\
	& IV        & 0.2      & 0.5  & 0.00978      & 0.00010     &  & -0.01915     & 0.00037    &  & -0.12970      & 0.01682     \\ \hline
	\end{tabular}
	\label{bayesa3}
\end{table}

\begin{table}[H]
	\centering
    \caption{{\scriptsize Bias and MSE of the Bayesian estimation by Metropolis-Hasting technique of $\hat{\beta}$ under squared error loss function $\hat{\beta}_{SEL}$ LINEX loss function $\hat{\beta}_{LL}$ and entropy loss function  $\hat{\beta}_{EL}$ for different choices of $n, m,\alpha,\beta$ and $T_1 = 0.4,T_2=4$.}}
\begin{tabular}{ccccrrrrrrrr}
\hline
& \multicolumn{2}{c}{} &      & \multicolumn{2}{c}{$\hat{\beta}_{SEL}$} &  & \multicolumn{2}{c}{$\hat{\beta}_{LL}$} &  & \multicolumn{2}{c}{$\hat{\beta}_{EL}$} \\ \cline{5-6} \cline{8-9} \cline{11-12} 
n,m       & sc   & $\alpha$    & $\beta$ & \multicolumn{1}{c}{Bias}         & \multicolumn{1}{c}{MSE}         &  & \multicolumn{1}{c}{Bias}         & \multicolumn{1}{c}{MSE}        &  & \multicolumn{1}{c}{Bias}          & \multicolumn{1}{c}{MSE}         \\ \hline
& I         & 0.2      & 0.5  & -0.02572     & 0.00066     &  & -0.04132     & 0.00171    &  & -0.03905      & 0.00153     \\
& II        & 0.2      & 0.5  & 0.06254      & 0.00391     &  & 0.04549      & 0.00207    &  & 0.04602       & 0.00212     \\
15,5  & III       & 0.2      & 0.5  & -0.02232     & 0.00050     &  & -0.05789     & 0.00335    &  & -0.13166      & 0.01734     \\
& IV        & 0.2      & 0.5  & -0.02203     & 0.00049     &  & -0.03712     & 0.00138    &  & -0.03318      & 0.00110     \\
& I         & 0.2      & 0.5  & 0.06556      & 0.00430     &  & 0.04950      & 0.00245    &  & 0.05266       & 0.00277     \\
& II        & 0.2      & 0.5  & 0.04908      & 0.00241     &  & 0.03299      & 0.00109    &  & 0.03581       & 0.00128     \\
20,10 & III       & 0.2      & 0.5  & 0.04123      & 0.00170     &  & 0.02461      & 0.00061    &  & 0.02561       & 0.00066     \\
& IV        & 0.2      & 0.5  & 0.05221      & 0.00273     &  & 0.03508      & 0.00123    &  & 0.03523       & 0.00124     \\
& I         & 0.2      & 0.5  & 0.27074      & 0.07330     &  & 0.20603      & 0.04245    &  & -0.01829      & 0.00033     \\
& II        & 0.2      & 0.5  & 0.20663      & 0.04269     &  & 0.14835      & 0.02201    &  & -0.26525      & 0.07036     \\
30,15 & III       & 0.2      & 0.5  & -0.00897     & 0.00008     &  & -0.06634     & 0.00440    &  & -0.33857      & 0.11463     \\
& IV        & 0.2      & 0.5  & 0.09338      & 0.00872     &  & 0.04952      & 0.00245    &  & -0.20184      & 0.04074     \\ \hline
	\end{tabular}
	\label{bayesb3}
\end{table}

\begin{table}[H]
	\centering
	\caption{{\scriptsize Bias and MSE of the Bayesian estimation by Metropolis-Hasting technique of $\hat{\alpha}$ under squared error loss function $\hat{\alpha}_{SEL}$, LINEX loss function $\hat{\alpha}_{LL}$, and entropy loss function $\hat{\alpha}_{EL}$, for different choices of $n, m,\alpha,\beta$ and $T_1 = 1,T_2=7$.}}
	\begin{tabular}{ccccrrrrrrrr}
		\hline
		& \multicolumn{2}{c}{} &      & \multicolumn{2}{c}{$\hat{\alpha}_{SEL}$} &  & \multicolumn{2}{c}{$\hat{\alpha}_{LL}$} &  & \multicolumn{2}{c}{$\hat{\alpha}_{EL}$} \\ \cline{5-6} \cline{8-9} \cline{11-12} 
		n,m       & sc    & $\alpha$    & $\beta$ & \multicolumn{1}{c}{Bias}         & \multicolumn{1}{c}{MSE}         &  & \multicolumn{1}{c}{Bias}         & \multicolumn{1}{c}{MSE}        &  & \multicolumn{1}{c}{Bias}          & \multicolumn{1}{c}{MSE}         \\ \hline
		& I         & 0.2      & 0.5  & 0.00518      & 0.00003     &  & -0.01761     & 0.00031    &  & -0.16221      & 0.02631     \\
		& II        & 0.2      & 0.5  & 0.05806      & 0.00337     &  & 0.04225      & 0.00178    &  & 0.02617       & 0.00068     \\
		15,5 & III       & 0.2      & 0.5  & -0.00693     & 0.00005     &  & -0.03054     & 0.00093    &  & -0.16305      & 0.02658     \\
		& IV        & 0.2      & 0.5  & 0.05205      & 0.00271     &  & 0.03126      & 0.00098    &  & -0.01381      & 0.00019     \\
		& I         & 0.2      & 0.5  & 0.05050      & 0.00255     &  & 0.02748      & 0.00075    &  & -0.05710      & 0.00326     \\
		& II        & 0.2      & 0.5  & -0.04886     & 0.00239     &  & -0.06539     & 0.00428    &  & -0.08724      & 0.00761     \\
		20,10 & III       & 0.2      & 0.5  & 0.00867      & 0.00008     &  & -0.02407     & 0.00058    &  & -0.07979      & 0.00637     \\
		& IV        & 0.2      & 0.5  & 0.02826      & 0.00080     &  & -0.01598     & 0.00026    &  & -0.12982      & 0.01685     \\
		& I         & 0.2      & 0.5  & 0.04949      & 0.00245     &  & 0.02357      & 0.00056    &  & -0.12603      & 0.01588     \\
		& II        & 0.2      & 0.5  & -0.02273     & 0.00052     &  & -0.03918     & 0.00154    &  & -0.06612      & 0.00437     \\
		30,15 & III       & 0.2      & 0.5  & 0.03766      & 0.00142     &  & 0.01636      & 0.00027    &  & -0.05739      & 0.00329     \\
		& IV        & 0.2      & 0.5  & 0.04199      & 0.00176     &  & 0.01822      & 0.00033    &  & -0.08413      & 0.00708     \\ \hline
	\end{tabular}
	\label{bayesa4}
\end{table}

\begin{table}[H]
	\centering
	 \caption{{\scriptsize Bias and MSE of the Bayesian estimation by Metropolis-Hasting technique of $\hat{\beta}$ under squared error loss function $\hat{\beta}_{SEL}$ LINEX loss function $\hat{\beta}_{LL}$ and entropy loss function  $\hat{\beta}_{EL}$ for different choices of $n, m,\alpha,\beta$ and $T_1 = 1,T_2=7$.}}
\begin{tabular}{ccccrrrrrrrr}
	\hline
	& \multicolumn{2}{c}{} &      & \multicolumn{2}{c}{$\hat{\beta}_{SEL}$} &  & \multicolumn{2}{c}{$\hat{\beta}_{LL}$} &  & \multicolumn{2}{c}{ $\hat{\beta}_{EL}$} \\ \cline{5-6} \cline{8-9} \cline{11-12} 
	n,m       & sc    & $\alpha$    & $\beta$ & \multicolumn{1}{c}{Bias}         & \multicolumn{1}{c}{MSE}         &  & \multicolumn{1}{c}{Bias}         & \multicolumn{1}{c}{MSE}        &  & \multicolumn{1}{c}{Bias}          & \multicolumn{1}{c}{MSE}         \\ \hline
	& I         & 0.2      & 0.5  & -0.01451     & 0.00021     &  & -0.04061     & 0.00165    &  & -0.13918      & 0.01937     \\
	& II        & 0.2      & 0.5  & -0.12429     & -0.12429    &  & -0.14812     & 0.02194    &  & -0.18092      & 0.03273     \\
	15,5  & III       & 0.2      & 0.5  & -0.20364     & 0.04147     &  & -0.23067     & 0.05321    &  & -0.29499      & 0.08702     \\
	& IV        & 0.2      & 0.5  & 0.04736      & 0.00224     &  & -0.00273     & 0.00001    &  & -0.18747      & 0.03514     \\
	& I         & 0.2      & 0.5  & 0.13877      & 0.01926     &  & 0.09580      & 0.00918    &  & -0.04724      & 0.00223     \\
	& II        & 0.2      & 0.5  & 0.05261      & 0.00277     &  & 0.03182      & 0.00101    &  & 0.02369       & 0.00056     \\
	20,10 & III       & 0.2      & 0.5  & 0.00837      & 0.00007     &  & -0.00689     & 0.00005    &  & -0.00274      & 0.00001     \\
	& IV        & 0.2      & 0.5  & 0.08137      & 0.00662     &  & 0.06326      & 0.00400    &  & 0.06173       & 0.00381     \\
	& I         & 0.2      & 0.5  & 0.04209      & 0.00177     &  & 0.02625      & 0.00069    &  & 0.02952       & 0.00087     \\
	& II        & 0.2      & 0.5  & 0.08297      & 0.00688     &  & 0.06655      & 0.00443    &  & 0.06909       & 0.00477     \\
	30,15 & III       & 0.2      & 0.5  & 0.04989      & 0.00249     &  & 0.03292      & 0.00108    &  & 0.03339       & 0.00111     \\
	& IV        & 0.2      & 0.5  & 0.03840      & 0.00147     &  & 0.02289      & 0.00052    &  & 0.02695       & 0.00073     \\ \hline
	\end{tabular}
	\label{bayesb4}
\end{table}

\section{Real Data Analysis}
In this section, a set of actual data set is given for simulation and illustration. The data set represents the times of failures and running times for samples of devices from an eld-tracking study of a larger system. The data set was studied by \cite{data1} \cite{data2} \cite{data3} The data set has 30 observations and it is given below(see \cite{data}  )

2.75, 0.13, 1.47, 0.23, 1.81, 0.30, 0.65, 0.10, 3.00, 1.73, 1.06, 3.00, 3.00, 2.12, 3.00, 3.00,
3.00, 0.02, 2.61, 2.93, 0.88, 2.47, 0.28, 1.43, 3.00, 0.23, 3.00, 0.80, 2.45, 2.66.

We use Anderson-Darling(AD) test and Kolmogorov-Smirnov(KS) test for checking the goodness of fit. The value for the AD test statistic is 1.3748 and the corresponding P-value is 0.2093. The value for the KS test statistic is 0.21649 and the corresponding P-value is 0.1201. Since the P-value in two test modes is high $(P>0.05)$, we cannot reject the null hypothesis. To summarize, we have sufficient evidence to hold that the chen distribution can provide a suitable fit for this data set.

Then we consider the estimation of two parameters under improved adaptive Type-II progressive censored data. We have used the original data set of data to generate random sets of improved adaptive Type-II progressive censored samples. We considered different combinations for $m$, $T_1$ and $T_2$.

The random samples obtained using different schemes and different combinations of $m$, $T_1$ and $T_2$ are given in Table \ref{scheme} 

We have used the original set of data to generate random sets of improved adaptive Type-II progressive censored samples in Table \ref{realdata}. We consider the two time threshold are $(T_1,T_2)=(0.4,4)$ and $(T_1,T_2)=(1,7)$. Based on the above situation, we have obtained the MLE, interval estimation and Bayesian estimation of $\alpha$ and $\beta$ are given in Table \ref{rmle1}-\ref{rbayesb4}.

According to these tables, we come to the following conclusions:

In most case, with increase of sample size $m$, Biases and MSE of MLE decrease. with the increase of time threshold, Biases and MSE of MLE and Bayes estimation increases. Estimates obtained from the uniform censoring scheme seems better than the results from uneven censoring schemes.

 We observe that the Bayes estimation under the squared error loss function, LINEX loss function, and entropy loss function shows minimum MSE than the MLE, for most cases of $n=30$, $m=5,15,20$.

\begin{table}[H]
	\centering
	\caption{{\scriptsize Improved adaptive Type-II progressive censoring schemes with $n$ = 30.}}
	\begin{tabular}{ccl}
		\hline
		m & censoring scheme  & censoring number
		\\ \hline
		~ & I & r=(0,0,0,0,25) \\ 
		~ & II & r=(10,3,4,0,8) \\ 
		5 & III & r=(5,3,10,5,2) \\ 
		~ & IV & r=(5,5,5,5,5) \\ 
		~ & I & r=(15,0,0,0,0,0,0,0,0,0,0,0,0,0,0) \\ 
		~ & II & r=(0,0,0,0,0,0,0,0,0,0,0,0,0,0,15) \\ 
		15 & III & r=(0,0,0,0,0,0,0,15,0,0,0,0,0,0,0) \\ 
		~ & IV & r=(1,1,1,1,1,1,1,1,1,1,1,1,1,1,1) \\ 
		~ & I & r=(1,1,1,1,1,1,1,0,0,0,0,0,0,0,0,0,0,1,1,1) \\ 
		~ & II & r=(1,1,0,0,0,0,0,0,0,0,0,0,0,0,0,0,0,0,0,8) \\ 
		20 & III & r=(1,0,1,0,0,0,0,0,0,0,0,0,0,0,0,0,0,1,0,7) \\ 
		~ & IV & r=(1,0,0,0,0,0,0,0,0,0,0,0,0,0,0,0,0,0,0,9) \\ \hline
	\end{tabular}
	\label{scheme}
\end{table}

\begin{sidewaystable}[htp]
	\centering
	\caption{{\scriptsize Improved adaptive Type-II progressive censored samples (CS) generated from the data under different censoring schemes}}
	\begin{tabular}{ccccl}
		\hline
		~ & ~ & {\scriptsize $T_1,T_2=(0.4,4)$} & {\scriptsize $T_1,T_2=(1,7)$} & ~  \\ 
		m & CS & $k_1$,$k_2$ & $k_1$,$k_2$ & Data\\ \hline
		~ & I & (5,5) & (5,5) & 0.02,0.10,0.13,0.23,0.23 \\ 
		5 & II & (0,5) & (0,5) & 1.06,1.81,2.66,2.75,3.00 \\ 
		~ & III & (1,5) & (2,5) & 0.28,0.88,2.75,3.00,3.00 \\ 
		~ & IV & (1,5) & (1,5) & 0.28,1.43,2.47,3.00,3.00 \\ 
		~ & I & (0,15) & (0,15) & 2.12,2.45,2.47,2.61,2.66,2.75,2.93,3.00,3.00,3.00,3.00,3.00,3.00,3.00,3.00 \\ 
		15 & II & (7,15) & (10,15) & 0.02,0.10,0.13,0.23,0.23,0.28,0.30,0.65,0.80,0.88,1.06,1.43,1.47,1.73,1.81 \\ 
		~ & III & (7,15) & (7,15) & 0.02,0.10,0.13,0.23,0.23,0.28,0.30,3.00,3.00,3.00,3.00,3.00,3.00,3.00,3.00 \\ 
		~ & IV & (4,15) & (5,15) & 0.10,0.23,0.28,0.28,0.88,1.43,1.73,2.12,2.47,2.66,2.93,3.00,3.00,3.00,3.00 \\ 
		~ & I & (3,20) & (5,20) & 0.10,0.23,0.28,0.65,0.88,1.43,1.73,1.81,2.12,2.45,2.47,2.61,2.66,2.75,2.93,3.00,3.00,3.00,3.00,3.00 \\ 
		20 & II & (5,20) & (8,20) & 0.10,0.23,0.23,0.28,0.30,0.65,0.80,0.88,1.06,1.43,1.47,1.73,1.81,2.12,2.45,2.47,2.61,2.66,2.75,3.00 \\ 
		~ & III & (5,20) & (8,20) & 0.10,0.13,0.23,0.28,0.30,0.65,0.80,0.88,1.06,1.43,1.47,1.73,1.81,2.12,2.45,2.47,2.61,2.75,2.93,3.00 \\ 
		~ & IV & (6,20) & (9,20) & 0.02,0.13,0.23,0.23,0.28,0.30,0.65,0.80,0.88,1.06,1.43,1.47,1.73,1.81,2.12,2.45,2.47,2.61,2.66,3.00 \\ \hline
	\end{tabular}
	\label{realdata}
\end{sidewaystable}


\begin{table}[H]
	\centering
	\caption{{\scriptsize Bias and MSE of the MLE $\hat{\alpha}$ and $\hat{\beta}$ for different $m$ of real data and $T_1=0.4,T_2=4$.}}
\begin{tabular}{ccccrrrrr}
	\hline
	& \multicolumn{2}{c}{} &      & \multicolumn{2}{c}{$\hat{\alpha}$} &  & \multicolumn{2}{c}{$\hat{\beta}$} \\ \cline{5-6} \cline{8-9} 
	m  & scheme    & $\alpha$    & $\beta$ & \multicolumn{1}{c}{Bias}         & \multicolumn{1}{c}{MSE}        &  & \multicolumn{1}{c}{Bias}        & \multicolumn{1}{c}{MSE}        \\ \hline
	& I         & 0.2      & 0.7  & -0.13755     & 0.01892    &  & -0.34460    & 0.11875    \\
	& II        & 0.2      & 0.7  & -0.16250     & 0.02640    &  & -0.08052    & 0.00648    \\
	5  & III       & 0.2      & 0.7  & -0.16375     & 0.02681    &  & -0.10764    & 0.01159    \\
	& IV        & 0.2      & 0.7  & -0.14705     & 0.02162    &  & -0.28734    & 0.08257    \\
	& I         & 0.2      & 0.7  & -0.11538     & 0.01331    &  & 0.02099     & 0.00044    \\
	& II        & 0.2      & 0.7  & 0.07061      & 0.00499    &  & -0.38778    & 0.15037    \\
	15 & III       & 0.2      & 0.7  & -0.12998     & 0.01689    &  & 0.06769     & 0.00458    \\
	& IV        & 0.2      & 0.7  & -0.14268     & 0.02036    &  & -0.22197    & 0.04927    \\
	& I         & 0.2      & 0.7  & 0.00829      & 0.00007    &  & -0.23010    & 0.05294    \\
	& II        & 0.2      & 0.7  & 0.01764      & 0.00031    &  & -0.18955    & 0.03593    \\
	20 & III       & 0.2      & 0.7  & -0.07151     & 0.00511    &  & 0.08936     & 0.00799    \\
	& IV        & 0.2      & 0.7  & -0.08887     & 0.00790    &  & 0.15267     & 0.02331    \\ \hline
	\end{tabular}
	\label{rmle1}
\end{table}

\begin{table}[H]
	\centering
	\caption{{\scriptsize  Bias and MSE of the MLE $\hat{\alpha}$ and $\hat{\beta}$ for different $m$ of real data and $T_1=1,T_2=7$.}}
\begin{tabular}{ccccrrrrr}
	\hline
	& \multicolumn{2}{c}{} &      & \multicolumn{2}{c}{$\hat{\alpha}$} &  & \multicolumn{2}{c}{$\hat{\beta}$} \\ \cline{5-6} \cline{8-9} 
	m  & scheme    & $\alpha$    & $\beta$ & \multicolumn{1}{c}{Bias}         & \multicolumn{1}{c}{MSE}        &  & \multicolumn{1}{c}{Bias}        & \multicolumn{1}{c}{MSE}        \\ \hline
	& I         & 0.2      & 0.7  & -0.14386     & 0.02070    &  & -0.38721    & 0.14993    \\
	& II        & 0.2      & 0.7  & -0.15644     & 0.02447    &  & -0.15489    & 0.02399    \\
	5  & III       & 0.2      & 0.7  & -0.18071     & 0.16391    &  & 0.03266     & 0.02687    \\
	& IV        & 0.2      & 0.7  & -0.16410     & 0.02693    &  & -0.10420    & 0.01086    \\
	& I         & 0.2      & 0.7  & -0.12755     & 0.01627    &  & 0.09112     & 0.00830    \\
	& II        & 0.2      & 0.7  & 0.06181      & 0.00382    &  & -0.32917    & 0.10835    \\
	15 & III       & 0.2      & 0.7  & -0.14848     & 0.02205    &  & -0.23134    & 0.05352    \\
	& IV        & 0.2      & 0.7  & -0.05637     & 0.00318    &  & -0.23491    & 0.05518    \\
	& I         & 0.2      & 0.7  & -0.02912     & 0.00085    &  & -0.14887    & 0.02216    \\
	& II        & 0.2      & 0.7  & -0.08256     & 0.00682    &  & 0.11204     & 0.01255    \\
	20 & III       & 0.2      & 0.7  & -0.08465     & 0.00717    &  & 0.11645     & 0.01356    \\
	& IV        & 0.2      & 0.7  & -0.09430     & 0.00889    &  & 0.16460     & 0.02709    \\ \hline
	\end{tabular}
	\label{rmle2}
\end{table}

\begin{table}[H]
	\centering
	\caption{{\scriptsize Interval estimations and Average Length(AL) for different $m$ of real data and $T_1=0.4,T_2=4.$}}
	\begin{tabular}{cccccc}
		\hline
		m & scheme & $\alpha$ & $\beta$ & $\alpha -AL$ & $\beta-AL$ \\ \hline
		~ & I & (0.01317,0.29599) & (0.10762,1.17368) & 0.28282  & 1.06606  \\ 
		~ & II & (0.01806,0.37642) & (0.03173,0.60056) & 0.35836  & 0.56882  \\ 
		5  & III & (0.00590,0.22284) & (0.21359,1.64287) & 0.21694  & 1.42928  \\ 
		~ & IV & (0.01260,0.33131) & (0.14081,1.20937) & 0.31870  & 1.06856  \\ 
		~ & I & (0.02120,0.49601) & (0.32968,1.57678) & 0.47482  & 1.24710  \\ 
		~ & II & (0.13796,0.53081) & (0.16223,0.72060) & 0.39284  & 0.55837  \\ 
		15  & III & (0.01777,0.27589) & (0.43800,1.34556) & 0.25812  & 0.90756  \\ 
		~ & IV & (0.03575,0.22436) & (0.26273,0.86978) & 0.18861  & 0.60705  \\ 
		~ & I & (0.10367,0.41849) & (0.28334,0.77930) & 0.31482  & 0.49596  \\ 
		~ & II & (0.11081,0.42743) & (0.31942,0.81573) & 0.31662  & 0.49631  \\ 
		20  & III & (0.05219,0.31635) & (0.53845,1.15719) & 0.26417  & 0.61874  \\ 
		~ & IV & (0.04259,0.28997) & (0.59634,1.21920) & 0.24738  & 0.62286 \\ \hline
	\end{tabular}
	\label{raci1}
\end{table}

\begin{table}[H]
	\centering
	\caption{{\scriptsize Interval estimations and Average Length(AL) for different $m$ of real data and $T_1=1,T_2=7.$}}
	\begin{tabular}{cccccc}
		\hline
		m & scheme & $\alpha$ & $\beta$ & $\alpha -AL$ & $\beta-AL$ \\ \hline
		~ & I & (0.01267,0.24871) & (0.09530,1.02666) & 0.23603  & 0.93137  \\ 
		~ & II & (0.00898,0.21137) & (0.19531,1.52142) & 0.20239  & 1.32611  \\ 
		5  & III & (0.00197,0.35623) & (0.38560,1.93552) & 0.35426  & 1.54991  \\ 
		~ & IV & (0.00585,0.36446) & (0.21602,1.64324) & 0.35861  & 1.42722  \\ 
		~ & I & (0.01482,0.55008) & (0.35590,1.75858) & 0.53526  & 1.40269  \\ 
		~ & II & (0.13254,0.51718) & (0.19312,0.71206) & 0.38464  & 0.51894  \\ 
		15  & III & (0.03205,0.20312) & (0.27755,0.79135) & 0.17107  & 0.51380  \\ 
		~ & IV & (0.06825,0.30224) & (0.26872,0.80495) & 0.23399  & 0.53623  \\ 
		~ & I & (0.08599,0.33957) & (0.34566,0.87874) & 0.25357  & 0.53308  \\ 
		~ & II & (0.04911,0.28086) & (0.56385,1.16947) & 0.23175  & 0.60562  \\ 
		20  & III & (0.04764,0.27928) & (0.61874,1.17494) & 0.23164  & 0.60760  \\ 
		~ & IV & (0.04104,0.27228) & (0.61047,1.22451) & 0.23125  & 0.61404 \\ \hline
	\end{tabular}
	\label{raci2}
\end{table}

\begin{table}[H]
	\centering
	\caption{{\scriptsize Bias and MSE of the Bayesian estimation by importance sampling technique of $\hat{\alpha}$ under squared error loss function $\hat{\alpha}_{SEL}$, LINEX loss function $\hat{\alpha}_{LL}$, and entropy loss function $\hat{\alpha}_{EL}$, for different $m$ of real data and $T_1 = 0.4,T_2=4$.}}
	\begin{tabular}{ccccrrrrrrrr}
		\hline
		& \multicolumn{2}{c}{} &      & \multicolumn{2}{c}{$\hat{\alpha}_{SEL}$} &  & \multicolumn{2}{c}{$\hat{\alpha}_{LL}$} &  & \multicolumn{2}{c}{ $\hat{\alpha}_{EL}$} \\ \cline{5-6} \cline{8-9} \cline{11-12} 
		m  & scheme    & $\alpha$    & $\beta$ & \multicolumn{1}{c}{Bias}       & \multicolumn{1}{c}{MSE}      &  & \multicolumn{1}{c}{Bias}      & \multicolumn{1}{c}{MSE}     &  & \multicolumn{1}{c}{Bias}       & \multicolumn{1}{c}{MSE}     \\ \hline
		& I         & 0.2      & 0.7  & -0.07636   & 0.00583  &  & -0.07787   & 0.00606  &  & -0.10805   & 0.01167  \\
		& II        & 0.2      & 0.7  & -0.08874   & 0.00787  &  & -0.09096   & 0.00827  &  & -0.11260   & 0.01268  \\
		5  & III       & 0.2      & 0.7  & -0.11380   & 0.01295  &  & -0.11393   & 0.01298  &  & -0.11576   & 0.01340  \\
		& IV        & 0.2      & 0.7  & -0.02445   & 0.00060  &  & -0.02474   & 0.00061  &  & -0.02715   & 0.00074  \\
		& I         & 0.2      & 0.7  & -0.06151   & 0.00378  &  & -0.06152   & 0.00378  &  & -0.06155   & 0.00379  \\
		& II        & 0.2      & 0.7  & -0.00017   & 0.00000  &  & -0.01691   & 0.00029  &  & -0.11025   & 0.01215  \\
		15 & III       & 0.2      & 0.7  & -0.05961   & 0.00355  &  & -0.06036   & 0.00364  &  & -0.06668   & 0.00445  \\
		& IV        & 0.2      & 0.7  & 0.01180    & 0.00014  &  & 0.01094    & 0.00012  &  & 0.00498    & 0.00002  \\
		& I         & 0.2      & 0.7  & 0.00445    & 0.00002  &  & 0.00424    & 0.00002  &  & 0.00264    & 0.00001  \\
		& II        & 0.2      & 0.7  & -0.02959   & 0.00088  &  & -0.03010   & 0.00091  &  & -0.03437   & 0.00118  \\
		20 & III       & 0.2      & 0.7  & 0.02763    & 0.00076  &  & 0.02709    & 0.00073  &  & 0.02399    & 0.00058  \\
		& IV        & 0.2      & 0.7  & 0.04536    & 0.00206  &  & 0.04257    & 0.00181  &  & 0.02717    & 0.00074  \\ \hline
	\end{tabular}
	\label{rbayesa1}
\end{table}

\begin{table}[H]
	\centering
    \caption{{\scriptsize Bias and MSE of the Bayesian estimation by importance sampling technique of $\hat{\beta}$ under squared error loss function $\hat{\beta}_{SEL}$ LINEX loss function $\hat{\beta}_{LL}$ and entropy loss function  $\hat{\beta}_{EL}$ for different $m$ of real data and $T_1 = 0.4,T_2=4$.}}
\begin{tabular}{ccccrrrrrrrr}
	\hline
	& \multicolumn{2}{c}{} &      & \multicolumn{2}{c}{$\hat{\beta}_{SEL}$} &  & \multicolumn{2}{c}{$\hat{\beta}_{LL}$} &  & \multicolumn{2}{c}{ $\hat{\beta}_{EL}$} \\ \cline{5-6} \cline{8-9} \cline{11-12} 
	m  & scheme    & $\alpha$    & $\beta$ & \multicolumn{1}{c}{Bias}        & \multicolumn{1}{c}{MSE}       &  &\multicolumn{1}{c}{Bias}        & \multicolumn{1}{c}{MSE}       &  & \multicolumn{1}{c}{Bias}        & \multicolumn{1}{c}{MSE}       \\ \hline
	& I         & 0.2      & 0.7  & -0.02299    & 0.00053    &  & -0.04440    & 0.00197    &  & -0.18193    & 0.03310    \\
	& II        & 0.2      & 0.7  & 0.03976     & 0.00158    &  & 0.03500     & 0.00122    &  & 0.02533     & 0.00064    \\
	5  & III       & 0.2      & 0.7  & 0.02922     & 0.00085    &  & 0.02915     & 0.00085    &  & 0.02906     & 0.00084    \\
	& IV        & 0.2      & 0.7  & 0.14193     & 0.02014    &  & 0.13866     & 0.01923    &  & 0.13424     & 0.01802    \\
	& I         & 0.2      & 0.7  & -0.01219    & 0.00015    &  & -0.01375    & 0.00019    &  & -0.01701    & 0.00029    \\
	& II        & 0.2      & 0.7  & 0.24726     & 0.06114    &  & 0.24455     & 0.05981    &  & 0.23998     & 0.05759    \\
	15 & III       & 0.2      & 0.7  & -0.05202    & 0.00271    &  & -0.05525    & 0.00305    &  & -0.06325    & 0.00400    \\
	& IV        & 0.2      & 0.7  & 0.19761     & 0.03905    &  & 0.19443     & 0.03780    &  & 0.18915     & 0.03578    \\
	& I         & 0.2      & 0.7  & 0.18049     & 0.03258    &  & 0.17883     & 0.03198    &  & 0.17615     & 0.03103    \\
	& II        & 0.2      & 0.7  & 0.02903     & 0.00084    &  & 0.02830     & 0.00080    &  & 0.02693     & 0.00073    \\
	20 & III       & 0.2      & 0.7  & -0.02497    & 0.00062    &  & -0.02515    & 0.00063    &  & -0.02551    & 0.00065    \\
	& IV        & 0.2      & 0.7  & 0.05447     & 0.00297    &  & 0.04743     & 0.00225    &  & 0.03455     & 0.00119    \\ \hline
	\end{tabular}
	\label{rbayesb1}
\end{table}

\begin{table}[H]
	\centering
		\caption{{\scriptsize Bias and MSE of the Bayesian estimation by importance sampling technique of $\hat{\alpha}$ under squared error loss function $\hat{\alpha}_{SEL}$, LINEX loss function $\hat{\alpha}_{LL}$, and entropy loss function $\hat{\alpha}_{EL}$, for different $m$ of real data and $T_1 = 1,T_2=7$.}}
	\begin{tabular}{ccccrrrrrrrr}
		\hline
		& \multicolumn{2}{c}{} &      & \multicolumn{2}{c}{$\hat{\alpha}_{SEL}$} &  & \multicolumn{2}{c}{$\hat{\alpha}_{LL}$} &  & \multicolumn{2}{c}{$\hat{\alpha}_{EL}$} \\ \cline{5-6} \cline{8-9} \cline{11-12} 
		m  & scheme    & $\alpha$    & $\beta$ & \multicolumn{1}{c}{Bias}        & \multicolumn{1}{c}{MSE}        &  & \multicolumn{1}{c}{Bias}        & \multicolumn{1}{c}{MSE}        &  & \multicolumn{1}{c}{Bias}         & \multicolumn{1}{c}{MSE}        \\ \hline
		& I         & 0.2      & 0.7  & -0.05041     & 0.00254    &  & -0.05144     & 0.00265    &  & -0.06200     & 0.00384    \\
		& II        & 0.2      & 0.7  & -0.10636     & 0.01131    &  & -0.10826     & 0.01172    &  & -0.12810     & 0.01641    \\
		5  & III       & 0.2      & 0.7  & 0.13604      & 0.01851    &  & 0.13582      & 0.01845    &  & 0.13487      & 0.01819    \\
		& IV        & 0.2      & 0.7  & 0.11202      & 0.01255    &  & 0.11174      & 0.01249    &  & 0.11044      & 0.01220    \\
		& I         & 0.2      & 0.7  & -0.05695     & 0.00324    &  & -0.05695     & 0.00324    &  & -0.05701     & 0.00325    \\
		& II        & 0.2      & 0.7  & 0.00840      & 0.00007    &  & -0.01088     & 0.00012    &  & -0.12055     & 0.01453    \\
		15 & III       & 0.2      & 0.7  & -0.02730     & 0.00075    &  & -0.02827     & 0.00080    &  & -0.03582     & 0.00128    \\
		& IV        & 0.2      & 0.7  & -0.01630     & 0.00027    &  & -0.01748     & 0.00031    &  & -0.02542     & 0.00065    \\
		& I         & 0.2      & 0.7  & -0.00692     & 0.00005    &  & -0.00730     & 0.00005    &  & -0.01038     & 0.00011    \\
		& II        & 0.2      & 0.7  & -0.06124     & 0.00375    &  & -0.06157     & 0.00379    &  & -0.06439     & 0.00415    \\
		20 & III       & 0.2      & 0.7  & 0.02467      & 0.00061    &  & -0.02508     & 0.00063    &  & -0.02842     & 0.00081    \\
		& IV        & 0.2      & 0.7  & 0.08591      & 0.00738    &  & 0.08308      & 0.00690    &  & 0.06954      & 0.00484    \\ \hline
	\end{tabular}
	\label{rbayesa2}
\end{table}

\begin{table}[H]
	\centering
    \caption{{\scriptsize Bias and MSE of the Bayesian estimation by importance sampling technique of $\hat{\beta}$ under squared error loss function $\hat{\beta}_{SEL}$ LINEX loss function $\hat{\beta}_{LL}$ and entropy loss function  $\hat{\beta}_{EL}$ for different $m$ of real data and $T_1 = 1,T_2=7$.}}
\begin{tabular}{ccccrrrrrrrr}
	\hline
	& \multicolumn{2}{c}{} &      & \multicolumn{2}{c}{$\hat{\beta}_{SEL}$} &  & \multicolumn{2}{c}{$\hat{\beta}_{LL}$} &  & \multicolumn{2}{c}{$\hat{\beta}_{EL}$} \\ \cline{5-6} \cline{8-9} \cline{11-12} 
	m  & scheme    & $\alpha$    & $\beta$ & \multicolumn{1}{c}{Bias}       & \multicolumn{1}{c}{MSE}        &  & \multicolumn{1}{c}{Bias}        & \multicolumn{1}{c}{MSE}        &  & \multicolumn{1}{c}{Bias}        & \multicolumn{1}{c}{MSE}        \\ \hline
	& I         & 0.2      & 0.7  & 0.05194     & 0.00270    &  & 0.03097     & 0.00096    &  & -0.15413    & 0.02376    \\
	& II        & 0.2      & 0.7  & 0.11752     & 0.01381    &  & 0.10943     & 0.01197    &  & 0.09549     & 0.00912    \\
	5  & III       & 0.2      & 0.7  & 0.00228     & 0.00001    &  & 0.00118     & 0.00000    &  & -0.00162    & 0.00000    \\
	& IV        & 0.2      & 0.7  & -0.07670    & 0.00588    &  & -0.07897    & 0.00624    &  & -0.08233    & 0.00678    \\
	& I         & 0.2      & 0.7  & 0.03894     & 0.00152    &  & 0.03823     & 0.00146    &  & 0.03671     & 0.00135    \\
	& II        & 0.2      & 0.7  & 0.22882     & 0.05236    &  & 0.22588     & 0.05102    &  & 0.22133     & 0.04899    \\
	15 & III       & 0.2      & 0.7  & 0.03251     & 0.00106    &  & 0.03111     & 0.00097    &  & 0.02666     & 0.00071    \\
	& IV        & 0.2      & 0.7  & 0.20779     & 0.04318    &  & 0.20631     & 0.04257    &  & 0.20332     & 0.04134    \\
	& I         & 0.2      & 0.7  & 0.01289     & 0.00017    &  & 0.01176     & 0.00014    &  & 0.00921     & 0.00008    \\
	& II        & 0.2      & 0.7  & -0.05215    & 0.00272    &  & -0.05243    & 0.00275    &  & -0.05302    & 0.00281    \\
	20 & III       & 0.2      & 0.7  & -0.04292    & 0.00184    &  & -0.04384    & 0.00192    &  & -0.04593    & 0.00211    \\
	& IV        & 0.2      & 0.7  & 0.01598     & 0.00026    &  & 0.01291     & 0.00017    &  & 0.00625     & 0.00004    \\ \hline
	\end{tabular}
	\label{rbayesb2}
\end{table}

\begin{table}[H]
	\centering
	\caption{{\scriptsize Bias and MSE of the Bayesian estimation by M-H technique of $\hat{\alpha}$ under squared error loss function $\hat{\alpha}_{SEL}$, LINEX loss function $\hat{\alpha}_{LL}$, and entropy loss function $\hat{\alpha}_{EL}$, for different $m$ of real data  and $T_1 = 0.4,T_2=4$.}}
\begin{tabular}{ccccrrrrrrrr}
	\hline
	& \multicolumn{2}{c}{} &      & \multicolumn{2}{c}{ $\hat{\alpha}_{SEL}$} &  & \multicolumn{2}{c}{$\hat{\alpha}_{LL}$} &  & \multicolumn{2}{c}{$\hat{\alpha}_{EL}$} \\ \cline{5-6} \cline{8-9} \cline{11-12} 
	m  & scheme    & $\alpha$    & $\beta$ & \multicolumn{1}{c}{Bias}         & \multicolumn{1}{c}{MSE}       &  & \multicolumn{1}{c}{Bias}         & \multicolumn{1}{c}{MSE}        &  & \multicolumn{1}{c}{Bias}         & \multicolumn{1}{c}{MSE}        \\ \hline
	& I         & 0.2      & 0.7  & -0.01205     & 0.00015    &  & -0.02702     & 0.00073    &  & -0.05351     & 0.00286    \\
	& II        & 0.2      & 0.7  & 0.05914      & 0.00350    &  & 0.01671      & 0.00028    &  & -0.07677     & 0.00589    \\
	5  & III       & 0.2      & 0.7  & 0.09361      & 0.00876    &  & 0.03098      & 0.00096    &  & -0.18152     & 0.03295    \\
	& IV        & 0.2      & 0.7  & 0.03262      & 0.00106    &  & -0.02590     & 0.00067    &  & -0.18804     & 0.03536    \\
	& I         & 0.2      & 0.7  & -0.01415     & 0.00020    &  & -0.02652     & 0.00070    &  & -0.02448     & 0.00060    \\
	& II        & 0.2      & 0.7  & 0.03592      & 0.00129    &  & -0.00947     & 0.00009    &  & -0.13514     & 0.01826    \\
	15 & III       & 0.2      & 0.7  & 0.08795      & 0.00773    &  & 0.06847      & 0.00469    &  & 0.03149      & 0.00099    \\
	& IV        & 0.2      & 0.7  & 0.11408      & 0.01301    &  & 0.09599      & 0.00921    &  & 0.07760      & 0.00602    \\
	& I         & 0.2      & 0.7  & 0.00804      & 0.00006    &  & -0.01747     & 0.00031    &  & -0.14197     & 0.02015    \\
	& II        & 0.2      & 0.7  & -0.00150     & 0.00000    &  & -0.02589     & 0.00067    &  & -0.12907     & 0.01666    \\
	20 & III       & 0.2      & 0.7  & -0.04960     & 0.00246    &  & -0.07071     & 0.00500    &  & -0.13727     & 0.01884    \\
	& IV        & 0.2      & 0.7  & -0.00149     & 0.00000    &  & -0.02496     & 0.00062    &  & -0.12762     & 0.01629    \\ \hline
	\end{tabular}
	\label{rbayesa3}
\end{table}

\begin{table}[H]
	\centering
	\caption{{\scriptsize Bias and MSE of the Bayesian estimates by M-H technique of $\hat{\beta}$ under squared error loss function $\hat{\beta}_{SEL}$ LINEX loss function $\hat{\beta}_{LL}$ and entropy loss function  $\hat{\beta}_{EL}$ for different $m$ of real data and $T_1 = 0.4,T_2=4$.}}
	\begin{tabular}{ccccrrrrrrrr}
		\hline
		& \multicolumn{2}{c}{} &      & \multicolumn{2}{c}{$\hat{\beta}_{SEL}$} &  & \multicolumn{2}{c}{$\hat{\beta}_{LL}$ } &  & \multicolumn{2}{c}{$\hat{\beta}_{EL}$} \\ \cline{5-6} \cline{8-9} \cline{11-12} 
		m  & scheme    & $\alpha$    & $\beta$ & \multicolumn{1}{c}{Bias}        & \multicolumn{1}{c}{MSE}        &  & \multicolumn{1}{c}{Bias}        & \multicolumn{1}{c}{MSE}        &  & \multicolumn{1}{c}{Bias}        & \multicolumn{1}{c}{MSE}        \\ \hline
		& I         & 0.2      & 0.7  & 0.00026     & 0.00000    &  & 0.16158     & 0.02611    &  & 0.11324     & 0.01282    \\
		& II        & 0.2      & 0.7  & 0.10786     & 0.01163    &  & 0.07988     & 0.00638    &  & 0.05756     & 0.00331    \\
		5  & III       & 0.2      & 0.7  & 0.08466     & 0.00717    &  & 0.05469     & 0.00299    &  & 0.03198     & 0.00102    \\
		& IV        & 0.2      & 0.7  & 0.17268     & 0.02982    &  & 0.14812     & 0.02194    &  & 0.13862     & 0.01922    \\
		& I         & 0.2      & 0.7  & -0.00555    & 0.00003    &  & -0.05769    & 0.00333    &  & -0.23209    & 0.05387    \\
		& II        & 0.2      & 0.7  & -0.00724    & 0.00005    &  & -0.04963    & 0.00246    &  & -0.14345    & 0.02058    \\
		15 & III       & 0.2      & 0.7  & 0.04112     & 0.00169    &  & -0.02387    & 0.00057    &  & -0.21480    & 0.04614    \\
		& IV        & 0.2      & 0.7  & 0.07185     & 0.00516    &  & 0.01602     & 0.00026    &  & -0.14434    & 0.02083    \\
		& I         & 0.2      & 0.7  & 0.08643     & 0.00747    &  & 0.02546     & 0.00065    &  & -0.17772    & 0.03158    \\
		& II        & 0.2      & 0.7  & 0.01869     & 0.00035    &  & -0.03623    & 0.00131    &  & -0.17325    & 0.03001    \\
		20 & III       & 0.2      & 0.7  & 0.05342     & 0.00285    &  & -0.00930    & 0.00009    &  & -0.16786    & 0.02818    \\
		& IV        & 0.2      & 0.7  & -0.01733    & 0.00030    &  & -0.08884    & 0.00789    &  & -0.29719    & 0.08832    \\ \hline
	\end{tabular}
	\label{rbayesb3}
\end{table}

\begin{table}[H]
	\centering
\caption{{\scriptsize Bias and MSE of the Bayesian estimates by M-H technique of $\hat{\alpha}$ under squared error loss function $\hat{\alpha}_{SEL}$, LINEX loss function $\hat{\alpha}_{LL}$, and entropy loss function $\hat{\alpha}_{EL}$, for different $m$ of real data  and $T_1 = 1,T_2=7$.}}
\begin{tabular}{ccccrrrrrrrr}
	\hline
	& \multicolumn{2}{c}{} &      & \multicolumn{2}{c}{$\hat{\alpha}_{SEL}$} &  & \multicolumn{2}{c}{$\hat{\alpha}_{LL}$} &  & \multicolumn{2}{c}{$\hat{\alpha}_{EL}$} \\ \cline{5-6} \cline{8-9} \cline{11-12} 
	m  & scheme    &$\alpha$    & $\beta$ & \multicolumn{1}{c}{Bias}         & \multicolumn{1}{c}{MSE}        &  & \multicolumn{1}{c}{Bias}         & \multicolumn{1}{c}{MSE}        &  & \multicolumn{1}{c}{Bias}         & \multicolumn{1}{c}{MSE}        \\ \hline
	& I         & 0.2      & 0.7  & -0.02627     & 0.00069    &  & -0.03927     & 0.00154    &  & -0.04494     & 0.00202    \\
	& II        & 0.2      & 0.7  & -0.03707     & 0.00137    &  & -0.06351     & 0.00403    &  & -0.11340     & 0.01286    \\
	5  & III       & 0.2      & 0.7  & 0.11327      & 0.01283    &  & 0.04831      & 0.00233    &  & -0.17754     & 0.03152    \\
	& IV        & 0.2      & 0.7  & 0.04528      & 0.00205    &  & -0.01225     & 0.00015    &  & -0.18208     & 0.03315    \\
	& I         & 0.2      & 0.7  & -0.00677     & 0.00005    &  & -0.01909     & 0.00036    &  & -0.01526     & 0.00023    \\
	& II        & 0.2      & 0.7  & 0.02574      & 0.00066    &  & -0.02209     & 0.00049    &  & -0.13430     & 0.01804    \\
	15 & III       & 0.2      & 0.7  & 0.11666      & 0.01361    &  & 0.09819      & 0.00964    &  & 0.07273      & 0.00529    \\
	& IV        & 0.2      & 0.7  & -0.02765     & 0.00076    &  & -0.05101     & 0.00260    &  & -0.16051     & 0.02577    \\
	& I         & 0.2      & 0.7  & 0.00706      & 0.00005    &  & -0.01802     & 0.00032    &  & -0.14993     & 0.02248    \\
	& II        & 0.2      & 0.7  & 0.05952      & 0.00354    &  & 0.03717      & 0.00138    &  & -0.06625     & 0.00439    \\
	20 & III       & 0.2      & 0.7  & 0.00884      & 0.00008    &  & -0.01686     & 0.00028    &  & -0.13304     & 0.01770    \\
	& IV        & 0.2      & 0.7  & 0.06622      & 0.00438    &  & 0.04601      & 0.00212    &  & -0.00263     & 0.00001    \\ \hline
	\end{tabular}
	\label{rbayesa4}
\end{table}

\begin{table}[H]
	\centering
		\caption{{\scriptsize  Bias and MSE of the Bayesian estimates by M-H technique of $\hat{\beta}$ under squared error loss function $\hat{\beta}_{SEL}$ LINEX loss function $\hat{\beta}_{LL}$ and entropy loss function  $\hat{\beta}_{EL}$ for different $m$ of real data and $T_1 = 1,T_2=7$.}}
	\begin{tabular}{ccccrrrrrrrr}
		\hline
		& \multicolumn{2}{c}{} &      & \multicolumn{2}{c}{$\hat{\beta}_{SEL}$} &  & \multicolumn{2}{c}{$\hat{\beta}_{LL}$} &  & \multicolumn{2}{c}{$\hat{\beta}_{EL}$} \\ \cline{5-6} \cline{8-9} \cline{11-12} 
		m  & scheme    & $\alpha$    & $\beta$ & \multicolumn{1}{c}{Bias}        & \multicolumn{1}{c}{MSE}        &  & \multicolumn{1}{c}{Bias}        & \multicolumn{1}{c}{MSE}        &  & \multicolumn{1}{c}{Bias}        & \multicolumn{1}{c}{MSE}        \\ \hline
		& I         & 0.2      & 0.7  & 0.08043     & 0.00647    &  & 0.25041     & 0.06271    &  & 0.22513     & 0.05068    \\
		& II        & 0.2      & 0.7  & 0.11075     & 0.01226    &  & 0.07756     & 0.00602    &  & 0.04768     & 0.00227    \\
		5  & III       & 0.2      & 0.7  & 0.03099     & 0.00096    &  & 0.00883     & 0.00008    &  & 0.00330     & 0.00001    \\
		& IV        & 0.2      & 0.7  & 0.11707     & 0.01371    &  & 0.08540     & 0.00729    &  & 0.05993     & 0.00359    \\
		& I         & 0.2      & 0.7  & 0.02181     & 0.00048    &  & -0.03164    & 0.00100    &  & -0.14292    & 0.02043    \\
		& II        & 0.2      & 0.7  & -0.03177    & 0.00101    &  & -0.07534    & 0.00568    &  & -0.16969    & 0.02879    \\
		15 & III       & 0.2      & 0.7  & 0.08504     & 0.00723    &  & 0.02962     & 0.00088    &  & -0.17040    & 0.02904    \\
		& IV        & 0.2      & 0.7  & 0.10244     & 0.01049    &  & 0.05754     & 0.00331    &  & -0.05893    & 0.00347    \\
		& I         & 0.2      & 0.7  & 0.09240     & 0.00854    &  & 0.03896     & 0.00152    &  & -0.08720    & 0.00760    \\
		& II        & 0.2      & 0.7  & 0.01618     & 0.00026    &  & -0.04503    & 0.00203    &  & -0.21788    & 0.04747    \\
		20 & III       & 0.2      & 0.7  & 0.00047     & 0.00000    &  & -0.05444    & 0.00296    &  & -0.22237    & 0.04945    \\
		& IV        & 0.2      & 0.7  & 0.05222     & 0.00273    &  & -0.00639    & 0.00004    &  & -0.15100    & 0.02280    \\ \hline
	\end{tabular}
	\label{rbayesb4}
\end{table}

\section{Conclusions}
In this paper, we have considered the problem of estimation of two unknown parameters $\alpha$ and $\beta$ based on IAT-II PCS sample from the Chen distribution, which has a bathtub-shape failure rate function. The maximum likelihood estimators and interval estimations by observed Fisher information matrix  of the two parameters have been obtained. On the premise that the prior distribution are gamma distribution, we mention that the Bayesian estimates are obtained under the squared error loss function, LINEX loss function, and entropy loss functions. To calculate the Bayesian estimates using the M-H technique and importance sampling technique. The results illustrate that in the simulation and real data tests, the maximum likelihood estimates and the Bayes estimate under different loss function  shown an insignificant difference, though the latter has slightly better performances than the former. Among the Bayes estimators of $\alpha$, estimator under squared loss function performs better and estimation effect under the entropy loss function comes second. Among the Bayes estimators of $\beta$, estimator under LINEX loss function possess minimum Bias and MSE.

In the future work, the parameter inference of IAT-II PCS under accelerated life test, application of IAT-II PCS in competitive failure test and parameter inference under more effective censoring scheme can be considered.

%
\section{Funding}
This research is supported by the National Natural Science Foundation of China.
\bibliographystyle{mystyle}
\bibliography{20220929}

\begin{thebibliography}{29}
\expandafter\ifx\csname natexlab\endcsname\relax\def\natexlab#1{#1}\fi
\providecommand{\url}[1]{\texttt{#1}}
\providecommand{\href}[2]{#2}
\providecommand{\path}[1]{#1}
\providecommand{\DOIprefix}{doi:}
\providecommand{\ArXivprefix}{arXiv:}
\providecommand{\URLprefix}{URL: }
\providecommand{\Pubmedprefix}{pmid:}
\providecommand{\doi}[1]{\href{http://dx.doi.org/#1}{\path{#1}}}
\providecommand{\Pubmed}[1]{\href{pmid:#1}{\path{#1}}}
\providecommand{\bibinfo}[2]{#2}
\ifx\xfnm\relax \def\xfnm[#1]{\unskip,\space#1}\fi
\bibitem[{$\ddot{o}$.G$\ddot{u}$r$\ddot{u}$nl$\ddot{u}$ Alma \&
  Belaghi(2016)}]{Alma2016}
\bibinfo{author}{$\ddot{o}$.G$\ddot{u}$r$\ddot{u}$nl$\ddot{u}$ Alma}, \&
  \bibinfo{author}{Belaghi, R.} (\bibinfo{year}{2016}).
\newblock \bibinfo{title}{On the estimation of the extreme value and normal
  distribution parameters based on progressive type-ii hybrid-censored data}.
\newblock {\it \bibinfo{journal}{Journal of Statistical Computation and
  Simulation}\/},  {\it \bibinfo{volume}{86}\/}, \bibinfo{pages}{569--596}.
\bibitem[{Balakrishnan \& Aggarwala(2000)}]{Balakrishnan2000}
\bibinfo{author}{Balakrishnan, N.}, \& \bibinfo{author}{Aggarwala, R.}
  (\bibinfo{year}{2000}).
\newblock \bibinfo{title}{Progressive censoring: Theory, methods, and
  applications}.
\bibitem[{Chen \& Bhattacharyya(1987)}]{Bhattacharyya1987}
\bibinfo{author}{Chen, S.-M.}, \& \bibinfo{author}{Bhattacharyya, G.~K.}
  (\bibinfo{year}{1987}).
\newblock \bibinfo{title}{Exact confidence bounds for an exponential parameter
  under hybrid censoring}.
\newblock {\it \bibinfo{journal}{Communications in Statistics - Theory and
  Methods}\/},  {\it \bibinfo{volume}{16}\/}, \bibinfo{pages}{2429--2442}.
\bibitem[{Chen(2000)}]{CHEN2000155}
\bibinfo{author}{Chen, Z.} (\bibinfo{year}{2000}).
\newblock \bibinfo{title}{A new two-parameter lifetime distribution with
  bathtub shape or increasing failure rate function}.
\newblock {\it \bibinfo{journal}{Statistics $\&$ Probability Letters}\/},  {\it
  \bibinfo{volume}{49}\/}, \bibinfo{pages}{155--161}.
\bibitem[{Childs et~al.(2008)Childs, Chandrasekar \&
  Balakrishnan}]{Childs2008ExactLI}
\bibinfo{author}{Childs, A.}, \bibinfo{author}{Chandrasekar, B.}, \&
  \bibinfo{author}{Balakrishnan, N.} (\bibinfo{year}{2008}).
\newblock \bibinfo{title}{Exact likelihood inference for an exponential
  parameter under progressive hybrid censoring schemes}.
\newblock (p. \bibinfo{pages}{323–334}).
\bibitem[{Childs et~al.(2003)Childs, Chandrasekar, Balakrishnan \&
  Kundu}]{childs2003}
\bibinfo{author}{Childs, A.}, \bibinfo{author}{Chandrasekar, B.},
  \bibinfo{author}{Balakrishnan, N.}, \& \bibinfo{author}{Kundu, D.}
  (\bibinfo{year}{2003}).
\newblock \bibinfo{title}{Exact likelihood inference based on type-i and
  type-ii hybrid censored samples from the exponential distribution}.
\newblock {\it \bibinfo{journal}{Annals of the Institute of Statistical
  Mathematics}\/}, .
\bibitem[{Cramer et~al.(2016)Cramer, Burkschat \& Górny}]{Cramer2016}
\bibinfo{author}{Cramer, E.}, \bibinfo{author}{Burkschat, M.}, \&
  \bibinfo{author}{Górny, J.} (\bibinfo{year}{2016}).
\newblock \bibinfo{title}{On the exact distribution of the mles based on
  type-ii progressively hybrid censored data from exponential distributions}.
\newblock {\it \bibinfo{journal}{Journal of Statistical Computation and
  Simulation}\/},  {\it \bibinfo{volume}{86}\/}, \bibinfo{pages}{2036--2052}.
\bibitem[{Dey \& Dey(2014)}]{dey2014}
\bibinfo{author}{Dey, S.}, \& \bibinfo{author}{Dey, T.} (\bibinfo{year}{2014}).
\newblock \bibinfo{title}{Statistical inference for the rayleigh distribution
  under progressively type-ii censoring with binomial removal}.
\newblock {\it \bibinfo{journal}{Applied Mathematical Modelling}\/},  {\it
  \bibinfo{volume}{38}\/}, \bibinfo{pages}{974--982}.
\bibitem[{El-Din et~al.(2017)El-Din, Amein, Shafay \& Mohamed}]{MohieElDin2017}
\bibinfo{author}{El-Din, M. M. M.~M.}, \bibinfo{author}{Amein, M.~M.},
  \bibinfo{author}{Shafay, A.~R.}, \& \bibinfo{author}{Mohamed, S.}
  (\bibinfo{year}{2017}).
\newblock \bibinfo{title}{Estimation of generalized exponential distribution
  based on an adaptive progressively type-ii censored sample}.
\newblock {\it \bibinfo{journal}{Journal of Statistical Computation and
  Simulation}\/},  {\it \bibinfo{volume}{87}\/}, \bibinfo{pages}{1292--1304}.
\bibitem[{Elshahhat \& Nassar(2022)}]{Ahmed2022}
\bibinfo{author}{Elshahhat, A.}, \& \bibinfo{author}{Nassar, M.}
  (\bibinfo{year}{2022}).
\newblock \bibinfo{title}{Analysis of adaptive type-ii progressively hybrid
  censoring with binomial removals}.
\newblock {\it \bibinfo{journal}{Journal of Statistical Computation and
  Simulation}\/},  {\it \bibinfo{volume}{0}\/}, \bibinfo{pages}{1--27}.
\bibitem[{Epstein(1954)}]{Epstein1954}
\bibinfo{author}{Epstein, B.} (\bibinfo{year}{1954}).
\newblock \bibinfo{title}{Truncated life tests in exponential case}.
\newblock {\it \bibinfo{journal}{The Annals of Mathematical Statistics}\/},
  {\it \bibinfo{volume}{25}\/}.
\bibitem[{Hastings(1970)}]{Hasting1970}
\bibinfo{author}{Hastings, W.~K.} (\bibinfo{year}{1970}).
\newblock \bibinfo{title}{{Monte Carlo sampling methods using Markov chains and
  their applications}}.
\newblock {\it \bibinfo{journal}{Biometrika}\/},  {\it \bibinfo{volume}{57}\/},
  \bibinfo{pages}{97--109}.
\bibitem[{Jerald.(1982)}]{1982}
\bibinfo{author}{Jerald., L.~F.} (\bibinfo{year}{1982}).
\newblock \bibinfo{title}{Statistical models and methods for lifetime data.}
\newblock {\it \bibinfo{journal}{John Wiley and Sons, New York}\/},  {\it
  \bibinfo{volume}{10}\/}, \bibinfo{pages}{316--317}.
\bibitem[{Kundu(2007)}]{Kundu2007}
\bibinfo{author}{Kundu, D.} (\bibinfo{year}{2007}).
\newblock \bibinfo{title}{On hybrid censored weibull distribution}.
\newblock {\it \bibinfo{journal}{Journal of Statistical Planning and
  Inference}\/},  {\it \bibinfo{volume}{137}\/}, \bibinfo{pages}{2127--2142}.
\bibitem[{Lv et~al.(2022)Lv, Tian \& Gui}]{qigui2022}
\bibinfo{author}{Lv, Q.}, \bibinfo{author}{Tian, Y.}, \& \bibinfo{author}{Gui,
  W.} (\bibinfo{year}{2022}).
\newblock \bibinfo{title}{Statistical inference for gompertz distribution under
  adaptive type-ii progressive hybrid censoring}.
\newblock {\it \bibinfo{journal}{Journal of Applied Statistics}\/},  {\it
  \bibinfo{volume}{0}\/}, \bibinfo{pages}{1--30}.
\bibitem[{Merovci \& Elbatal(2014)}]{data2}
\bibinfo{author}{Merovci, F.}, \& \bibinfo{author}{Elbatal, I.}
  (\bibinfo{year}{2014}).
\newblock \bibinfo{title}{Weibull-rayleigh distribution: Theory and
  applications}.
\newblock {\it \bibinfo{journal}{Applied Mathematics $\&$ Information
  Sciences}\/},  {\it \bibinfo{volume}{5}\/}.
\bibitem[{Metropolis et~al.(1953)Metropolis, Rosenbluth, Rosenbluth, Teller \&
  Teller}]{Metropolis1953}
\bibinfo{author}{Metropolis, N.}, \bibinfo{author}{Rosenbluth, A.~W.},
  \bibinfo{author}{Rosenbluth, M.~N.}, \bibinfo{author}{Teller, A.~H.}, \&
  \bibinfo{author}{Teller, E.} (\bibinfo{year}{1953}).
\newblock \bibinfo{title}{Equation of state calculations by fast computing
  machines}.
\newblock {\it \bibinfo{journal}{Journal of Chemical Physics}\/},  {\it
  \bibinfo{volume}{21}\/}, \bibinfo{pages}{1087--1092}.
\bibitem[{Mohan \& Chacko(2020)}]{data3}
\bibinfo{author}{Mohan, R.}, \& \bibinfo{author}{Chacko, M.}
  (\bibinfo{year}{2020}).
\newblock \bibinfo{title}{Estimation of parameters of kumaraswamy-exponential
  distribution based on adaptive type-ii progressive censored schemes}.
\newblock {\it \bibinfo{journal}{Journal of Statistical Computation and
  Simulation}\/},  {\it \bibinfo{volume}{91}\/}, \bibinfo{pages}{81 -- 107}.
\bibitem[{Mohan \& Chacko(2021)}]{Rakhi2021}
\bibinfo{author}{Mohan, R.}, \& \bibinfo{author}{Chacko, M.}
  (\bibinfo{year}{2021}).
\newblock \bibinfo{title}{Estimation of parameters of kumaraswamy-exponential
  distribution based on adaptive type-ii progressive censored schemes}.
\newblock {\it \bibinfo{journal}{Journal of Statistical Computation and
  Simulation}\/},  {\it \bibinfo{volume}{91}\/}, \bibinfo{pages}{81--107}.
\bibitem[{Nassar \& Abo-Kasem(2017)}]{Nassar2017}
\bibinfo{author}{Nassar, M.}, \& \bibinfo{author}{Abo-Kasem, O.}
  (\bibinfo{year}{2017}).
\newblock \bibinfo{title}{Estimation of the inverse weibull parameters under
  adaptive type-ii progressive hybrid censoring scheme}.
\newblock {\it \bibinfo{journal}{Journal of Computational and Applied
  Mathematics}\/},  {\it \bibinfo{volume}{315}\/}, \bibinfo{pages}{228--239}.
\bibitem[{Nelson(1998)}]{data1}
\bibinfo{author}{Nelson, W.} (\bibinfo{year}{1998}).
\newblock \bibinfo{title}{Statistical methods for reliability data. $
  \mathrm{New}$ $\mathrm{York(NY)}$: Wiley}.
\bibitem[{Ng et~al.(2009)Ng, Kundu \& Chan}]{Ng2009}
\bibinfo{author}{Ng, H. K.~T.}, \bibinfo{author}{Kundu, D.}, \&
  \bibinfo{author}{Chan, P.~S.} (\bibinfo{year}{2009}).
\newblock \bibinfo{title}{Statistical analysis of exponential lifetimes under
  an adaptive type-ii progressive censoring scheme}.
\newblock {\it \bibinfo{journal}{Naval Research Logistics (NRL)}\/},  {\it
  \bibinfo{volume}{56}\/}, \bibinfo{pages}{687--698}.
\bibitem[{Oguntunde et~al.(2017)Oguntunde, Adejumo \& Owoloko}]{data}
\bibinfo{author}{Oguntunde, P.~E.}, \bibinfo{author}{Adejumo, A.~O.}, \&
  \bibinfo{author}{Owoloko, E.~A.} (\bibinfo{year}{2017}).
\newblock \bibinfo{title}{Application of kumaraswamy inverse exponential
  distribution to real lifetime data}.
\newblock {\it \bibinfo{journal}{International journal of applied mathematics
  and statistics}\/},  {\it \bibinfo{volume}{56}\/}, \bibinfo{pages}{34--47}.
\bibitem[{Panahi \& Moradi(2020)}]{Hanieh2020}
\bibinfo{author}{Panahi, H.}, \& \bibinfo{author}{Moradi, N.}
  (\bibinfo{year}{2020}).
\newblock \bibinfo{title}{Estimation of the inverted exponentiated rayleigh
  distribution based on adaptive type ii progressive hybrid censored sample}.
\newblock {\it \bibinfo{journal}{Journal of Computational and Applied
  Mathematics}\/},  {\it \bibinfo{volume}{364}\/}, \bibinfo{pages}{112345}.
\bibitem[{Tomer \& Panwar(2015)}]{Sanjeev2015}
\bibinfo{author}{Tomer, S.~K.}, \& \bibinfo{author}{Panwar, M.~S.}
  (\bibinfo{year}{2015}).
\newblock \bibinfo{title}{Estimation procedures for maxwell distribution under
  type-i progressive hybrid censoring scheme}.
\newblock {\it \bibinfo{journal}{Journal of Statistical Computation and
  Simulation}\/},  {\it \bibinfo{volume}{85}\/}, \bibinfo{pages}{339--356}.
\bibitem[{Yan et~al.(2021)Yan, Li \& Yu}]{YAN202138}
\bibinfo{author}{Yan, W.}, \bibinfo{author}{Li, P.}, \& \bibinfo{author}{Yu,
  Y.} (\bibinfo{year}{2021}).
\newblock \bibinfo{title}{Statistical inference for the reliability of burr-xii
  distribution under improved adaptive type-ii progressive censoring}.
\newblock {\it \bibinfo{journal}{Applied Mathematical Modelling}\/},  {\it
  \bibinfo{volume}{95}\/}, \bibinfo{pages}{38--52}.
\bibitem[{Ye et~al.(2014)Ye, Chan, Xie \& Ng}]{Ye2014}
\bibinfo{author}{Ye, Z.-S.}, \bibinfo{author}{Chan, P.-S.},
  \bibinfo{author}{Xie, M.}, \& \bibinfo{author}{Ng, H. K.~T.}
  (\bibinfo{year}{2014}).
\newblock \bibinfo{title}{Statistical inference for the extreme value
  distribution under adaptive type-ii progressive censoring schemes}.
\newblock {\it \bibinfo{journal}{Journal of Statistical Computation and
  Simulation}\/},  {\it \bibinfo{volume}{84}\/}, \bibinfo{pages}{1099--1114}.
\bibitem[{Zhang \& Gui(2020)}]{Zhanggui2020}
\bibinfo{author}{Zhang, Y.}, \& \bibinfo{author}{Gui, W.}
  (\bibinfo{year}{2020}).
\newblock \bibinfo{title}{A goodness of fit test for the pareto distribution
  with progressively type ii censored data based on the cumulative hazard
  function}.
\newblock {\it \bibinfo{journal}{Journal of Computational and Applied
  Mathematics}\/},  {\it \bibinfo{volume}{368}\/}, \bibinfo{pages}{112557}.
\bibitem[{Zhang \& Gui(2019)}]{Zhangandgui2019}
\bibinfo{author}{Zhang, Z.}, \& \bibinfo{author}{Gui, W.}
  (\bibinfo{year}{2019}).
\newblock \bibinfo{title}{Statistical inference of reliability of generalized
  rayleigh distribution under progressively type-ii censoring}.
\newblock {\it \bibinfo{journal}{Journal of Computational and Applied
  Mathematics}\/},  {\it \bibinfo{volume}{361}\/}, \bibinfo{pages}{295--312}.

\end{thebibliography}

\end{document}